\documentclass[11pt]{article}

 \usepackage{benstyle}
 
\usepackage{appendix}
 
 \usepackage{enumerate}
\usepackage{longtable}

\usepackage{multirow,float}

\DefineNamedColor{named}{Purple}{cmyk}{0.45,0.86,0,0}

\newcommand*{\per}{\mathrm{per}}
\newcommand*{\fold}{\mathrm{fold}}
\newcommand*{\ivl}{\mathrm{int}}
\newcommand*{\lt}{\mathrm{left}}
\newcommand*{\rt}{\mathrm{right}}
\newcommand*{\supp}{\mathrm{supp}}

\title{On stable reconstructions from nonuniform Fourier measurements}
\author{Ben Adcock \\ Department of Mathematics \\ Purdue University \\ \hspace{20pc} \\ \vspace{-2.25pc} \\
150 N. University Street \\ West Lafayette, IN 47907 \\ USA  \and Milana Gataric \\  CCA, Centre for Mathematical Sciences \\ \hspace{20pc} \\ \vspace{-2.25pc} \\ University of Cambridge \\ Wilberforce Rd, Cambridge CB3 0WA \\ United Kingdom \and Anders C. Hansen \\ DAMTP, Centre for Mathematical Sciences \\ University of Cambridge \\ Wilberforce Rd, Cambridge CB3 0WA \\ United Kingdom }

\begin{document}

\maketitle

\begin{abstract}
We consider the problem of recovering a compactly-supported function from a finite collection of pointwise samples of its Fourier transform taking nonuniformly.  First, we show that under suitable conditions on the sampling frequencies -- specifically, their density and bandwidth -- it is possible to recover any such function $f$ in a stable and accurate manner in any given finite-dimensional subspace; in particular, one which is well suited for approximating $f$.  In practice, this is carried out using so-called nonuniform generalized sampling (NUGS).  Second, we consider approximation spaces in one dimension consisting of compactly supported wavelets.  We prove that a linear scaling of the dimension of the space with the sampling bandwidth is both necessary and sufficient for stable and accurate recovery.  Thus wavelets are up to constant factors optimal spaces for reconstruction.
\end{abstract}

\section{Introduction}\label{s:introduction}
The reconstruction of an image or signal from its Fourier measurements is an important task in applied mathematics and engineering.  It arises in numerous applications, ranging from Magnetic Resonance Imaging (MRI) to X-ray Computed Tomography (CT), seismology and microscopy (the latter processes usually involve the Radon transform, but equate to reconstruction from Fourier measurements via the Fourier slice theorem).

The purpose of this paper is to consider the following question: given fixed measurements of an unknown image $f$, i.e.\ a finite collection of samples of its Fourier transform $\hat{f}$, not necessarily taken on a Cartesian grid, under what conditions is possible to recover an approximation to $f$ in a given finite-dimensional space $\rT$, and how can this be achieved with a stable numerical algorithm?  Our main contributions are: (i) a theoretical framework for understanding when stable reconstruction is possible, along with a stable numerical algorithm to achieve such a reconstruction, (ii) a full answer to the question of stable recovery in the univariate case based on the density and bandwidth of the samples, and (iii) analysis for the applications-relevant case where the reconstruction space $\rT$ consists of wavelets.

The particular focus of this paper is on the case where the data is acquired nonuniformly in the Fourier domain.  Nonuniform sampling arises naturally in many of the applications listed above.  In particular, radial sampling of the Fourier transform results whenever sampling with the Radon transform.  Furthermore, nonuniform sampling patterns -- in particular, spiral trajectories (see \cite{AhnEtAlSpiralNMR,DelattreEtAlSprial,PruessmannUnserMRIFast,StanfordMRI,MeyerEtAlSpiralCoronary} and references therein) -- have become increasingly popular in MRI applications in the last several decades. 

\subsection{Generalized sampling}

The approach we take in this paper is based on recent work in sampling and reconstruction in abstract Hilbert spaces, known as \textit{generalized sampling} (GS).  GS, in the form we consider in this paper, was introduced by two of the authors in \cite{BAACHShannon}.  Yet, its roots can be traced to earlier work of Unser \& Aldroubi \cite{unser1994general}, Eldar \cite{eldar2003sampling}, Eldar \& Werther \cite{eldar2005general}, Gr\"ochenig \cite{GrochenigIrregular,GrochenigIrregularExpType,GrochenigModernSamplingBook}, Shizgal \& Jung \cite{shizgalGegen2}, Hrycak \& Gr\"ochenig \cite{hrycakIPRM} and others.  See also the work of Aldroubi and others on average sampling \cite{AldroubiAverageSamp,AldroubiConvolutionAvg,SunAvgFRI,SunXhouAverages}.  GS addresses the following problem in sampling theory.  Suppose that a finite number of samples of an element $f$ of a Hilbert space are given as inner products with respect to a particular basis or frame.  Suppose also that $f$ can be efficiently represented in another basis or frame (e.g.\ it has sparse or rapidly-decaying coefficients).  GS obtains a reconstruction of $f$ in this new system using only the original data.  In the linear case, this is achieved by least-squares fitting \cite{BAACHAccRecov}, but when sparsity is assumed, one can combine it with compressed sensing techniques (i.e.\ convex optimization and random sampling) to achieve substantial subsampling \cite{BAACHGSCS}.  By doing so, one obtains techniques for infinite-dimensional (i.e.\ analog) compressed sensing, known as GS--CS.

GS/GS--CS has been considered for the Fourier reconstruction problem whenever the samples are taken uniformly \cite{BAACHShannon,BAACHAccRecov}.  The primary advantage of this approach over more standard reconstruction algorithms in medical imaging (e.g.\ gridding or iterative reconstructions -- see below) is that it allows one to take advantage of the availability of efficient representation systems for images.  It is well known that natural images are well represented using wavelets.  Images may be sparse in wavelets, or have coefficients with rapid decay.  Moreover, representing medical images in such systems has other benefits over classical Fourier series representations, such as improved compressibility, better feature detection and easier and more effective denoising \cite{LaineWaveletsBiomed,NowakWaveletDenoise,WeaverEtAlWaveletFiltering}.  GS allows one to compute quasi-optimal reconstructions in wavelets from the given set of Fourier samples, and therefore exploit such beneficial properties.  In the case of uniform Fourier samples, the use of GS/GS--CS with wavelets was studied in \cite{AHPWavelet,AHPRBreaking}.

\subsection{Contributions of the paper}
The focus of this paper is the case of GS with nonuniform Fourier measurements.  Although a particular instance of GS corresponding to (weighted) Fourier frames, we shall refer to the resulting framework as \textit{nonuniform generalized sampling} (NUGS) for the purposes of expositional clarity. 

Specifically, suppose that $\Omega = \{ \omega_1,\ldots,\omega_N \} \allowbreak  \subseteq \bbR^d$ is a set of $N$ frequencies in $d \geq 1$ dimensions, and that we are given the measurements $\{ \hat{f}(\omega) : \omega \in \Omega \}$ of an unknown signal $f \in \rL^2(D)$, where $D \subseteq \bbR^d$ is compact.  Note that for the purposes of this paper, $\Omega$ is fixed and cannot readily be altered (this is typical in the applications listed above).  Let $\rT \subseteq \rL^2(D)$ be a finite-dimensional space in which we wish to recover $f$.  For example, $\rT$ could consist of the first $M$ functions in some wavelet basis.   In this paper, we derive conditions on $\Omega$ and $\rT$ under which stable reconstruction is possible with NUGS.  In the one-dimensional setting we show that if the samples $\Omega$ have \textit{density} $\delta < 1$ then stable reconstruction is possible provided the \textit{bandwidth} $K$ of $\Omega$ is sufficiently large, with the precise nature of this scaling depending on the properties of $\rT$.  We also address the case of the critical density $\delta = 1$ within the context of Fourier frames.

An important facet of NUGS is that it is always possible to compute the various constants that enter into the stability and error estimates.  Thus, for a given $\Omega$ and $\rT$, stable reconstruction can be guaranteed \textit{a priori} by a straightforward numerical calculation.  Our numerical results illustrate that these bounds give good estimates of the actual reconstruction errors seen in practice.

In \S \ref{s:recon_prob} and \S \ref{s:NUGS} we present general theory for NUGS, and in  \S \ref{s:samp_thm} we give the classification of stable recovery in the univariate setting in terms of bandwidth and density.  In the second half of the paper, \S \ref{s:wavelets}--\S \ref{s:num_exp}, we address the case where the subspace $\rT$ corresponds to a wavelet basis.  A result proved in \cite{AHPWavelet} ($d=1$) and \cite{AHKM2DWavelets} ($d=2$) shows that when the sampling set $\Omega$ consists of the first $N$ uniform frequencies one can recover the first $\ord{N}$ coefficients in an arbitrary wavelet basis with GS.  Thus wavelet bases are, up to constants, optimal bases in which to recover images from uniform Fourier samples.  This is not true for example with algebraic polynomial bases, in which case one can stably recover only the first $\ordu{\sqrt{N}}$ coefficients (see \cite{hrycakIPRM}, as well as \cite{BAACHAccRecov,AdcockHansenShadrinStabilityFourier}).  In \S \ref{s:wavelets} we extend the $d=1$ result to the nonuniform case.  Specifically, if the samples $\Omega$ have density $\delta < 1$ and bandwidth $K>0$ then we prove that one can recover the first $\ord{K}$ wavelet coefficients stably and accurately.  Thus there is a one-to-one relationship between the sampling bandwidth and the wavelet scale.  This is further highlighted in \S \ref{s:bandwidth}, where we show that any attempt to reconstruct a fixed number of wavelet coefficients from a sampling bandwidth $K$ below a critical threshold (that is linear in $K$) necessarily results in exponential ill-conditioning.  This generalizes a result first proved in \cite{AHPWavelet} for uniform samples.

Let us now make several further remarks.  In \S \ref{s:samp_thm}--\S\ref{s:num_exp} we focus on one-dimensional functions.  Whilst our main motivations arise from two- or three-dimensional imaging problems -- in particular, MRI and X-ray CT -- we make this simplification so that the paper is of reasonable length and self-contained.  Generalizations of the main results proved in \S \ref{s:samp_thm}--\S\ref{s:num_exp} are currently in progress, and will be reported on in the near future (see \S \ref{s:conclusions} for a discussion).  Having said this, we note that there are also one-dimensional applications in imaging, in which case when one's interest lies only in particular cross-sectional slices of the whole image.  For example, see the work in Nuclear Magnetic Resonance (NMR) spectroscopy \cite{MRI1D2,  NMR3, NMRbook,  NMR2, NMR1} and also in MRI \cite{MRI1D1, MRI1D3}. 
One-dimensional approaches for acquiring wavelet coefficients also find use in MRI, such as in wavelet-encoding techniques \cite{HealyWeaverWaveletIEEE,LaineWaveletsBiomed,PanychWaveletEncoding,PanychEtAlWaveletEncoded,WeaverEtAlWaveletFiltering,WeaverEtAlWaveletEncoding}.  Hence a one-dimensional study is not just for reasons of brevity, but has also practical relevance. 

We remark also that this paper does not address the issue of sparsity.  Sparsity-exploiting algorithms are currently revolutionizing signal and image reconstruction.  Since a main focus of this paper is wavelets, in which images are known to be sparse, it may at first sight appear strange not to seek to exploit such properties.  For uniform samples this has indeed been done by using the aforementioned GS--CS framework, and the results are reported in \cite{BAACHGSCS,AHPRBreaking}.  However, as was explained in \cite{BAACHGSCS} (see also \cite{AHRT2013AEIP}), before one can exploit sparsity it is first necessary to understand the underlying linear mapping between the samples and coefficients in the reconstruction system.  This is precisely what we do in this paper via NUGS.  Exploiting sparsity by extending the work of \cite{AHPRBreaking} to the case of fully nonuniform Fourier samples is a topic of ongoing investigations.

 \subsection{Relation to previous work}\label{ss:relation}
Two well-known algorithms in MRI reconstruction are gridding \cite{JacksonEtAlGridding,SedaratDCFOptimal,GelbNonuniformFourier} and resampling \cite{RosenfeldURS1}.  Our work differs from both in that we assume an analog model for the image $f$, as opposed to viewing $f$ as a finite-length Fourier series.  Consequently, a key issue in NUGS is that of \textit{approximation}.  By using an appropriate reconstruction space $\rT$, we avoid the unpleasant artefacts (e.g.\ ringing) associated with these algorithms.

Another popular method for MRI reconstruction is the iterative reconstruction algorithm \cite{KnoppKunisPotts,MatejFesslerIterative,FesslerFastIterativeMRI}.  As we explain, this is a special case of our framework corresponding to a pixel basis for $\rT$.  Thus our work provides as a corollary theoretical guarantees for the stability and error of this algorithm. Moreover, by changing $\rT$ (which is permissible given a fixed number of samples due to the result proved in \S \ref{s:wavelets}), we can exploit the advantages of higher-order wavelets to obtain better reconstructions.

This work is also related to a large body of previous research into nouniform sampling theory.  Often in nonuniform sampling, one models the samples as giving rise to a Fourier frame \cite{AldroubiGrochenigSIREV,BenedettoIrregular,BenedettoSpiral,FeichtingerGrochenigIrregular}.  Reconstruction can then be carried out by iterative inversion of the frame operator, for example (this reconstruction is quite different from ours, though, since the approximation properties are tied to those of the sampling frame).  This can be problematic for several reasons.  First, even if a sequence of samples gives rise to a frame, it can be difficult to determine the frame constants so as to get explicit bounds.  Moreover, in our setting, where we consider finite sets of samples, and reconstructions in finite-dimensional spaces, the existence of a countable frame sequence is in some senses superfluous.  Instead, we shall mainly focus on simple conditions for stable recovery -- namely, the density $\delta$ and the bandwidth $K$ -- both of which can be easily computed.

Second, Fourier frames do not allow arbitrary clustering of sampling frequencies, such as often the case in practice.  To address this, the usual approach in the literature, which we shall also use in this paper, is to use ``adaptive weights'' \cite{FeichtingerGrochenigIrregular,FeichtingerEtAlEfficientNonuniform,GrochenigIrregular,GrochenigLAA,GrochenigIrregularExpType,GrochenigModernSamplingBook,GrochenigStrohmerMarvasti} (also referred to as ``density compensation factors'').  See also \cite{KunisPottsScattered} and references therein.  In \cite{FeichtingerGrochenigIrregular,GrochenigModernSamplingBook,GrochenigStrohmerMarvasti} an efficient algorithm  for the nonuniform sampling problem was introduced, known as the ACT algorithm (Adaptive weights, Conjugate gradients, Toeplitz).  In the terminology of this paper, it reconstructs in 
a subspace $\rT$ consisting of shifted Diracs, which corresponds to a trigonometric polynomial approximation of $\hat{f}$.  Subject to appropriate density conditions, in \cite{GrochenigIrregularExpType,GrochenigModernSamplingBook} Gr\"ochenig proves convergence  of the approximation to $f$ as the number of samples increases.  Note that the focus of this work is slightly different.  Gr\"ochenig et al.\ primarily consider the recovery of a bandlimited function from nonuniform pointwise samples, whereas we consider the recovery of a compactly-supported function from pointwise samples of its Fourier transform.  Although mathematically equivalent, the setup affects the choice of reconstruction space.  In our setting for example, a Dirac basis would not be ideal for approximating an image $f$, whereas wavelet bases are typically well suited.  Having said this, our paper can also be seen as an extension of this work, where the primary innovations are the generalization to arbitrary subspaces $\rT$ and the detailed analysis for the specific case of wavelet spaces.

\section{The reconstruction problem}\label{s:recon_prob}
We first introduce some notation.  Throughout the paper, we write $\ip{\cdot}{\cdot}$ for the inner product on $\rL^2(\bbR^d)$ and $\nm{\cdot}$ for the corresponding norm.  We denote the Fourier transform by 
\bes{
\hat{f}(\omega) = \int_{\bbR^d} f(x) \E^{-2 \pi \I \omega \cdot x} \D x,\quad \omega \in \bbR^d,\quad f \in \rL^2(\bbR^d).
}
Our primary concern are functions with compact support in a domain $D \subseteq \bbR^d$.  Thus, we define the subspace 
$
\rH = \left \{ f \in \rL^2(\bbR^d) :\ \mathrm{supp}(f) \subseteq D \right \} \subseteq \rL^2(\bbR^d).
$
We shall use the notation $\Omega = \{ \omega_1,\ldots,\omega_N \} \subseteq \bbR^d$
to denote a finite set of distinct frequencies, henceforth referred to as a \textit{sampling scheme}.  We also define $\rT$ to be a finite-dimensional subspace of $\rH$; the so-called \textit{reconstruction space}.  The corresponding orthogonal projection onto $\rT$ is denoted by $\cP_{\rT}$.

Given a sampling scheme $\Omega$ and reconstruction space $\rT \subseteq \rH$, we wish to compute an approximation $\tilde{f} \in \rT$ to $f$ in the subspace $\rT$ using only the sampling data 
\be{\label{smaplingdata}
\{ \hat{f}(\omega) : \omega \in \Omega  \}.
}
As discussed in \cite{BAACHOptimality}, when developing a method for this problem, i.e.\ a mapping $F:f\mapsto\tilde{f}$ depending solely on the data \R{smaplingdata}, there are two critical considerations:
\begin{enumerate}[(i)]
\item The mapping $F$ should be \textit{quasi-optimal}: for some constant $\mu = \mu(F) \ll \infty$, we have $\| f - F(f) \| \leq \mu \| f - \cP_{\rT} f \|$, $\forall f \in \rH$.  Recall that the motivation for considering a particular reconstruction space $\rT$ is that $f$ is known to be well-represented in this space.  In other words, the error $\| f - \cP_{\rT} f \|$ is small.  Quasi-optimality guarantees that the reconstruction $\tilde{f}$ from the data \R{smaplingdata} inherits such a small error.  
\item The mapping $F$ should be \textit{numerically stable}: for some constant $\kappa = \kappa(F) \ll \infty$, we have $\| F(g) \| \leq \kappa \| g \|$, $\forall g \in \rH$.  This property is of course vital to ensure that perturbations of the measurements do not adversely affect the reconstruction.
\end{enumerate}
With this to hand, the main focus of the paper is to answer the following questions: (i) under what conditions on $\Omega$ and $\rT$ stable, quasi-optimal reconstruction is possible, and (ii) how large is the reconstruction constant $C = C(\Omega,\rT)$.  We do this by analyzing a particular instance of GS, so-called NUGS, which we introduce in \S \ref{s:NUGS}.  This provides a sufficient condition for (i) and an upper bound for (ii). As we explain in Remark \ref{r:optimality}, however, under appropriate conditions the NUGS reconstruction cannot be outperformed by any other method.  Hence our analysis of NUGS provides not only sufficient conditions for stable, quasi-optimal reconstruction, but also (under appropriate, but mild, assumptions) necessary conditions.

\section{Generalized sampling for nonuniform Fourier samples}\label{s:NUGS}

Suppose $\Omega = \{ \omega_1,\ldots,\omega_N \}$ is a sampling scheme.  In the case of uniform Fourier samples taken at the Nyquist rate, stability and accuracy of the GS reconstruction is guaranteed by strong convergence of the sampling operator $\cS : \rH \rightarrow \rH$ defined by
\be{
\label{FS}
\cS f(x) = \sum^{N}_{n=1} \hat{f}(\omega_n) \E^{2 \pi \I \omega_n \cdot x} \bbI_{D}(x),
}
to the identity operator as $N \rightarrow \infty$ \cite{BAACHOptimality}.  When considering nonuniform samples, we shall use the following weaker condition:

\defn{
\label{d:admissible}
Let $\Omega$ be a sampling scheme, $\cS :  \rH\rightarrow \rH$ a bounded linear operator and let $\rT$ be a finite-dimensional subspace of $\rH$. Suppose that $\cS$ satisfies
\begin{enumerate}[I]
\item for each $f \in \rH$, $\cS f$ depends only on the sampling data $\{ \hat{f}(\omega) : \omega \in \Omega \}$,
\item $\cS$ is self-adjoint with respect to $\ip{\cdot}{\cdot}$ and satisfies
\be{
\label{triangle}
|\ip{\cS f}{g}| \leq \sqrt{\ip{\cS f}{f} \ip{\cS g}{g} },\quad \forall f,g \in \rH,
}
\item there exists a positive constant $C_1 = C_1(\Omega,\rT)$ such that
\be{
\label{approxParseval1}
\ip{\cS f}{f} \geq C_1 \| f \|^2,\quad \forall f \in \rT.
}
\end{enumerate}
Then $\cS$ is said to be an admissible sampling operator for the pair $(\Omega,\rT)$.
}

For convenience, throughout the remainder of the paper we shall assume that $C_1$ is the largest constant for which \R{approxParseval1} holds.  Given such an operator $\cS$, we now also define the constants $C_2 = C_2(\Omega)$ and $C_{3} = C_3(\Omega,\rT)$ by
\ea{
\label{approxParseval2}
\ip{\cS f}{f} &\leq C_2 \| f \|^2,\quad \forall f \in \rH,
\\
\label{approxParseval3}
\ip{\cS f}{f} &\leq C_3 \| f \|^2,\quad \forall f \in \rT.
}
Likewise, we assume these constants are the smallest possible.  Note that $C_2$ and $C_3$ exist since $\cS$ is bounded, and we also trivially have that $C_3 \leq C_2$.  Typically, we want $C_2$ to be independent of the number of samples $N$ (or more precisely, the bandwidth $K$ of $\Omega$ -- see Definition \ref{d:K_delta}), since, as we see later, it appears in the error and stability estimates.

Note that the inequalities \R{approxParseval1} and \R{approxParseval3} ensure that the bilinear form $\ip{\cS \cdot}{\cdot}$ gives rise to an equivalent inner product on $\rT$.  This relaxes the condition of strong convergence of the operator $\cS$.  We remark in passing that in the case of uniform sampling, the operator $\cS$ defined by \R{FS} is automatically an admissible sampling operator whenever $N$ is sufficiently large, with constants $C_1 \approx 1$ for large $N$ and $C_2 = 1$ for all $N$.

In the nonuniform setting, there are many potential ways to construct the operator $\cS$.  In this paper, we focus primarily on the following simple construction:
\be{
\label{weightedFS}
\cS f(x) = \sum^{N}_{n=1} \mu_{n} \hat{f}(\omega_{n}) \E^{2 \pi \I \omega_n \cdot x} \bbI_{D}(x),
}
where $\mu_{n} > 0$ are particular weights.  Observe that $\cS$, when defined in this way, automatically satisfies properties (i) and (ii) for an admissible sampling operator.  Clearly, in the case of uniform sampling, \R{weightedFS} reduces to \R{FS} when the weights $\mu_n =1$.

Given a sampling scheme $\Omega$, a finite-dimensional subspace $\rT$ and an admissible sampling operator $\cS$ we now define the GS reconstruction by
\be{
\label{GSnonuniform}
\tilde f \in \rT,\qquad \ip{\cS \tilde f}{g} = \ip{\cS f}{g},\quad \forall g \in \rT,
}
and write $F = F_{\Omega,\rT}$ for the mapping $f \mapsto \tilde{f}$.  Note that if $\cS$ is given by \R{weightedFS} then this is equivalent to the weighted least-squares data fit:
\be{
\label{GS_LS_data}
\tilde{f} =  \underset{g \in \rT}{\operatorname{argmin}} \sum^{N}_{n=1} \mu_n \left | \hat{f}(\omega_n) - \hat{g}(\omega_n) \right |^2.
}
Although \R{GSnonuniform} is an instance of GS corresponding to nonuniform Fourier samples, we shall refer to it as \textit{nonuniform generalized sampling (NUGS)} for the purposes of clarity.  

As we shall see next, the constants $C_1$ and $C_2$ arising from an admissible sampling operator $\cS$ determine the stability and quasi-optimality of the resulting NUGS reconstruction.  We first define the corresponding reconstruction constant $C(\Omega,\rT)$:

\defn{\label{defn_constant}
Let $\cS$ be an admissible sampling operator with constants $C_1$ and $C_2$ given by $\R{approxParseval1}$ and $\R{approxParseval2}$ respectively.  The ratio $C(\Omega,\rT) = \sqrt{C_2 / C_1}$ is referred to as the NUGS reconstruction constant.
}

\thm{
\label{t:bounds}
Let $\Omega$ be a sampling scheme and $\rT$ a finite-dimensional subspace, and suppose that $\cS$ is an admissible sampling operator.  Then the reconstruction $F(f) = \tilde{f}$ defined by $\R{GSnonuniform}$ exists uniquely for any $f \in \rH$ and we have the sharp bound
\be{
\label{strong_fNbd}
\| f - F(f+h) \| \leq \tilde{C} \left ( \| f - \cP_{\rT} f \| + \| h \| \right ),\quad \forall f,h \in \rH,
}
where the constant $\tilde{C}$ is given by
$
\tilde{C} = \tilde{C}(\Omega,\rT) =  \left \{ \| g \|/\| \cP_{\cS(\rT)} g \| : g \in \rT, g \neq 0 \right \}.
$
Moreover, the constant $\tilde{C}$ satisfies $\tilde{C} \leq C$, where $C = C(\Omega,\rT)$ is the corresponding reconstruction constant $($Definition $\ref{defn_constant}$$)$.
}

This theorem is a particular instance for the case of nonuniform samples of a result proved for GS in \cite{BAACHOptimality} (for completeness, we include a simplified proof in the appendix).  It confirms that admissibility of $\cS$ is sufficient for quasi-optimality and stability of the reconstruction $\tilde{f}$ up to the magnitude of the reconstruction constant $C$.   Note that the result is true under the slightly weaker assumption $\tilde{C} < \infty$ (which is of course implied by $C_1>0$ and $C_2<\infty$).  However, the constant $\tilde{C}$ is rather difficult to work with in practice \cite{BAACHOptimality}. 

\rem{ \label{r:S_choice}
Although we assume throughout the remainder of the paper that $\cS$ takes the form $\R{weightedFS}$, the results of this section do not require this.  They only assume that $\cS$ is admissible in the sense of Definition $\ref{d:admissible}$.  This allows one to consider more general forms for $\cS$ than the diagonal choice $\R{weightedFS}$, as has recently been considered in several works.  In $\textnormal{\cite{GelbSongFrameNFFT}}$, Gelb \& Song use banded operators $\cS$ for nonuniform Fourier sampling, and in $\textnormal{\cite{BergerGrochenigOblique}}$ Berger \& Gr\"ochenig consider improved choices for $\cS$ within the setting of GS in general Hilbert spaces.
}

\subsection{Computation of the reconstruction}\label{ss:computation}
We now discuss implementation of NUGS.  Since this is similar to the the more general case of GS, we give only a brief summary (see \cite{BAACHOptimality} for details).  Recall that if $\cS$ is given by \R{weightedFS}, then \R{GSnonuniform} is equivalent to \R{GS_LS_data}.  In particular, if $\{ \phi_m\}^{P}_{m=1}$ is a basis for $\rT$, and if the reconstruction $\tilde f$ is given by $f = \sum^{P}_{m=1} a_m \phi_m$, then the vector $a = (a_1,\ldots,a_P)^{\top}$ is the least squares solution of the $N \times P$ linear system
$
A a \approx b,
$
where $b=(b_1,\ldots,b_N)^{\top}$ and $A \in \bbC^{N \times P}$ have entries
\be{
\label{Abdef}
b_n = \sqrt{\mu_{n}}  \hat{f}(\omega_{n}),\quad A_{n,m} = \sqrt{\mu_{n}} \widehat{\phi_m}(\omega_{n}),\quad n=1,\ldots,N,\  m=1,\ldots,P.
}
Hence, once a basis for $\rT$ is specified, $\tilde f$ can be computed by solving a least squares problem.  The computational cost in doing so is proportional to the condition number $\kappa(A)$, which determines the number of iterations required in an iterative solver such as conjugate gradients,  multiplied by the cost of performing matrix-vector operations with $A$ and its adjoint $A^*$.  Note that if an orthonormal basis is specified for $\rT$, then $\kappa(A) \leq C(\Omega,\rT)$ \cite{BAACHOptimality}.

\rem{
\label{r:wavelet_fast}
Efficient computation of $\tilde{f}$ relies on a fast algorithm for  performing efficient matrix-vector computations. The existence of such algorithms, however, depends solely on the choice of the reconstruction space $\rT$. In the case of wavelet reconstruction spaces, $\ord{N \log N}$ algorithms can be constructed, based on Nonuniform Fast Fourier Transforms $($NUFFTs$)$ $\textnormal{\cite{PottsEtAlNUFFTtutorial,Potts}}$.  See $\S \ref{s:num_exp}$ for further discussion.

Recall from $\S \ref{ss:relation}$ that the ACT algorithm $\textnormal{\cite{FeichtingerEtAlEfficientNonuniform,GrochenigIrregularExpType,GrochenigModernSamplingBook}}$ can be viewed as an instance of NUGS where $\hat{\rT} = \{ \hat{g} : g \in \rT \}$ is a space of trigonometric polynomials on a compact interval.  Efficient implementation in $\ord{N \log N}$ time is carried out using fast Toeplitz solvers, although one could also use NUFFTs with the same overall complexity (see $\textnormal{\cite{KnoppKunisPotts}}$), as we shall do in the case of wavelet choices for $\rT$.
}

\section{A generalized sampling theorem for univariate nonuniform samples}\label{s:samp_thm}

In this section, we provide a generalized sampling theorem which asserts that stable, quasi-optimal reconstruction is possible for any fixed $\rT$ under appropriate conditions on the nonuniform sampling scheme $\Omega$.  Note that our focus will now be on the case $d =1$ -- see \S \ref{s:conclusions} for further discussion on the extension to $d \geq 1$.  We shall consider two scenarios in  the next two subsections.  First, sampling schemes $\Omega$ subject to appropriate density and bandwidth conditions.  Second, sampling schemes arising from Fourier frames.  

\subsection{$(K,\delta)$-dense sampling schemes}\label{ss:dense_stable}
We commence with the following definition:
\defn{
\label{d:K_delta}
Let $K >0$, $0< \delta < 1$ and $\omega_1 < \omega_2 < \ldots < \omega_N$.  The sampling scheme $\Omega = \{ \omega_1,\ldots,\omega_N \}$ has bandwidth $K$ and density $\delta$ if $\Omega \subseteq [-K,K]$ and
\bes{
\max_{n=0,\ldots,N} \left \{ \omega_{n+1} - \omega_n \right \} \leq \delta,
}
where $\omega_0 = \omega_N - 2 K$ and $\omega_{N+1} = \omega_{1} + 2K$.  In this case, we say that $\Omega$ is $(K,\delta)$-dense.
}

Our main result in this section is to show that, for an arbitrary fixed reconstruction space $\rT$, $(K,\delta)$-density for suitably large $K$ and small $\delta$ ensures stable reconstruction.  This holds provided the weights $\mu_n$ in \R{weightedFS} are chosen according to the following strategy:
\be{
\label{weights}
\mu_n = \tfrac12 \left ( \omega_{n+1} - \omega_{n-1} \right ),\quad n=1,\ldots,N,
}
where, as above, we set $\omega_0 = \omega_N  - 2 K$ and $\omega_{N+1} = \omega_{1} + 2 K$ (other choices of weights are possible, but we shall not address this issue). Note that our $(K, \delta)$-density condition is similar to the definition of $\delta$-density used by Gr\"ochenig \cite{GrochenigIrregular,GrochenigLAA} and Feichtinger, Gr\"{o}chenig \& Strohmer \cite{FeichtingerEtAlEfficientNonuniform}, the only difference being we make the bandwidth $K$ explicit. Also, note that the weights \R{weights} are exactly the same as those used therein.

We now require the following lemma:

\lem{
\label{l:FT_decay}
Let $\Omega = \{ \omega_1,\ldots,\omega_N \}$ be $(K,\delta)$-dense and suppose that $\mu_1,\ldots,\mu_N$ are given by $\R{weights}$.  Then for any nonzero $f \in \rL^2(0,1)$ we have
\bes{
 \left ( \sqrt{1-\| \hat{f} \|^2_{\bbR \backslash I} /\| f \|^2 } - \delta \right )^2  \| f \|^2 \leq \sum^{N}_{n=1} \mu_n | \hat{f}(\omega_n) |^2 \leq (1+\delta)^2 \| f \|^2,
}
where $I = (-K + \delta/2 , K - \delta/2 )$, and $\nmu{\hat{f}}_{\bbR \backslash I}^{2} = \int_{\bbR \backslash I} | \hat{f}(\omega) |^2 \D \omega$. 
}

This lemma is an extension -- with a similar proof -- of a result of Gr\"ochenig \cite{GrochenigIrregular} to the case where the number of samples $N$ is finite.  Gr\"ochenig's result is obtained in the limit $N,K \rightarrow \infty$.  We also note that the lower bound is strictly less than $(1-\delta)^2$ for any nonzero $f$, since $f$ is compactly supported and hence $\hat{f}$ cannot have compact support.  However, the lower bound converges to $(1-\delta)^2$ as the bandwidth $K$ is increased.  In other words, $N$ Fourier samples with density $\delta < 1$ and appropriately large bandwidth $K$ are sufficient to control $\| f \|$.  This observation will lead to the main result in this section.

\prf{
Define the function $F\in\rL^2\left(-1/2,1/2\right)$ by $F(x) = f(x+1/2)$. Since $| \hat{F}(\omega) | = | \hat{f}(\omega) |$, and also $\|F\|=\|f\|$, it is enough to prove the theorem for $F$.  Let $z_{n} = \frac{1}{2}(\omega_{n-1} + \omega_n)$ and write $\chi(\omega) = \sum^{N}_{n=1} \hat{F}(\omega_n) \bbI_{[z_n,z_{n+1})}(\omega)$ so that
\bes{
S^2 = \sum^{N}_{n=1} \mu_n | \hat{F}(\omega_n) |^2 = \int^{z_{N+1}}_{z_{1}} | \chi(x) |^2 \D x = \| \chi \|_{J}^{2},
}
where $J=(z_1,z_{N+1})$ and $\nm{\cdot}_{J}$ denotes the $L^2$-norm over $J$.  Hence
\be{
\label{S_ineq}
\| \hat{F} \|_{J} - \| \hat{F} - \chi \|_{J}   \leq S \leq \| \hat{F} \|_{\bbR} + \| \hat{F} - \chi \|_{J}.
}
Using Wirtinger's inequality \cite[Lem.\ 1]{GrochenigIrregular}, we find that
\eas{
\| \hat{F} - \chi \|^2_{J} &= \sum^{N}_{n=1} \int^{z_{n+1}}_{z_n} \left | \hat{F}(\omega) - \hat{F}(\omega_n) \right |^2 \D \omega
\\
& =  \sum^{N}_{n=1} \left ( \int^{\omega_n}_{z_{n}} + \int^{z_{n+1}}_{\omega_n} \right ) \left | \hat{F}(\omega) - \hat{F}(\omega_n) \right |^2 \D \omega
\\
& \leq \sum^{N}_{n=1} \left ( \frac{4(\omega_n - z_n)^2}{\pi^2} \int^{\omega_n}_{z_n}+ \frac{4(z_{n+1}-\omega_n)^2}{\pi^2} \int^{z_{n+1}}_{\omega_n} \right )  \left | \frac{\D }{\D \omega} \hat{F}(\omega) \right |^2  \D \omega\\
& \leq \frac{\delta^2}{\pi^2} \int_{J} \left | \frac{\D }{\D \omega} \hat{F}(\omega) \right |^2 \D \omega,
}
where the final inequality follows from the $(K,\delta)$-density of the samples.  Since differentiation in Fourier space corresponds to multiplication by $(-2 \pi \I  x)$ in physical space, we conclude that
$
\nmu{ \hat{F} - \chi }_{J} \leq 2\delta \nmu{ \widehat{F_1 } }_{J} \leq 2\delta \nmu{ \widehat{F_1 } }_{\bbR},
$
where  $F_1(x) = x F(x).$
Since $F$ is supported in $[-1/2,1/2]$, we deduce that
\be{
\label{err_1}
\| \hat{F}-\chi \|_{J} \leq 2\delta \|F_1 \| \leq \delta \| F \|.
}
Substituting this into the right-hand side of \R{S_ineq} gives $S \leq (1+\delta) \| F \|$, and hence the upper bound.  For the lower bound, we first note that $I \subseteq J$.  Hence, by \R{S_ineq} and \R{err_1},
\eas{
S \geq \| \hat{F} \|_{I} - \delta \| F \| \geq \sqrt{\| \hat{F} \|^2 - \| \hat{F} \|^2_{\bbR \backslash I} } -\delta \| F \|,
}
and the lower bound follows.
}

\defn{
Let $\rT \subseteq \rH$.  The $z$-residual of $\rT$ is the quantity
\be{
\label{z_residual}
E(\rT,z) = \sup \left \{ \| \hat{f} \|_{\bbR \backslash (-z,z)} : f \in \rT,\  \| f \| =1 \right \} ,\quad z \in [0,\infty).
}
} 
Note that $E(\rT,z) \leq 1$, $\forall z$ and any $\rT$, since $\| \hat{f} \| = \| f \|$.

\lem{
Let $\rT \subseteq \rH$ be a finite-dimensional subspace.  Then $E(\rT,z) \rightarrow 0$ monotonically as $z \rightarrow \infty$.
}
\prf{
Clearly $E(\rT,z)$ is monotonically decreasing in $z$.  Moreover, for any fixed $f \in \rT$, we have $\| \hat{f} \|_{\bbR \backslash (-z,z)} \rightarrow 0$ as $z \rightarrow \infty$.  The result now follows immediately.
}

Combining the previous two lemmas, we immediately obtain our main result of this section:

\thm{
\label{t:Cbound_arbitrary}
Let $\rT \subseteq \rH$ be finite-dimensional and let $\Omega$ be $(K,\delta)$-dense, where
\bes{
\delta < \sqrt{1-E(\rT,K-1/2)^2}.  
}
Let $\cS$ be given by $\R{weightedFS}$ with weights $\R{weights}$.  Then $\cS$ is admissible with reconstruction constant $C(\Omega,\rT)$ (see Definition $\ref{defn_constant}$) satisfying
\be{
\label{Cs_bound}
C(\Omega,\rT) \leq \frac{1+\delta}{\sqrt{1-E(\rT,K-1/2)^2 } - \delta}.
}
}

\prf{
The upper bound in Lemma \ref{l:FT_decay} immediately gives $C_2(\Omega) \leq (1+\delta)^2$.  For $C_1(\Omega,\rT)$ we set $f = g \in \rT$ in Lemma \ref{l:FT_decay}, and then apply the definition of $E(\rT,z)$ to get
\bes{
C_{1}(\Omega,\rT) \geq \left ( \sqrt{1-E(\rT,K-\delta/2)^2} - \delta \right )^2.
}
The result now follows from monotonicity of $E(\rT,z)$ and the definition of $C(\Omega,\rT)$.
}

This theorem states the following.  For a fixed reconstruction space $\rT$, the reconstruction constant $C(\Omega,\rT)$ can be made arbitrarily close to $\frac{1+\delta}{1-\delta}$ by taking $K$ sufficiently large.  Thus, even with highly nonuniform samples, we are guaranteed a stable reconstruction for large enough bandwidth $K$ provided the density condition $\delta < 1$ holds, with the precise level of stability controlled primarily by how close $\delta$ is to one.  As noted previously, in \cite{GrochenigIrregular} it was shown that infinite sequences $\{ \omega_n \}_{n \in \bbN}$ with bandwidth $K = \infty$ and density $\delta < 1$ are sets of sampling, i.e.\ they give rise to weighted Fourier frames $\{ \sqrt{\mu_n} \E^{2 \pi \I \omega_n \cdot} \bbI_{[0,1]}(\cdot) \}_{n \in \bbN}$ for $\rH$.  Therefore, based on arguments given in  \cite{GrochenigIrregular}, Theorem \ref{t:Cbound_arbitrary} shows that this condition also allows one to stably reconstruct from finitely-many samples in any finite-dimensional subspace $\rT$, provided the sampling bandwidth is sufficiently large.

A key aspect of the Theorem \ref{t:Cbound_arbitrary} is the nature of the bound \R{Cs_bound}.  The right-hand side separates geometric properties of the sampling scheme $\Omega$, i.e.\ the density $\delta$, from intrinsic properties of the reconstruction space $\rT$, i.e.\ the $z$-residual $E(\rT,z)$.  Hence, by analyzing the $z$-residual for each particular choice of $\rT$, we can guarantee stable, quasi-optimal reconstruction for \textit{all} sampling schemes $\Omega$ with $\delta <1$ and appropriate bandwidth $K$.  This is how we shall proceed in \S \ref{s:wavelets} when we provide recovery guarantees for wavelet reconstruction spaces.  We note in passing that a universal lower bound for $E(\rT,z)$ for any subspace $\rT$ of dimension $M$ is provided by the $M^{\rth}$ eigenvalue of the prolate spheroidal wavefunctions \cite{ProlatesIII}.  In particular, ensuring $E(\rT,z) < c$ for some $c<1$ necessitates at least a linear scaling of $z$ with $M$, regardless of the choice of $\rT$.  For wavelets, we show that a linear scaling is also sufficient.  

\rem{
In $\textnormal{\cite{GrochenigIrregularExpType}}$, Gr\"ochenig proves stability and convergence of the aforementioned ACT algorithm.  As mentioned, this algorithm can be seen as a particular case corresponding to a trigonometric basis in frequency.  The contribution of Theorem $\ref{t:Cbound_arbitrary}$ is that it allows for arbitrary spaces $\rT$.  Note that in Gr\"ochenig's case (up to some minor differences in how the boundary is dealt with), $E(\rT,K-1/2) = 0$ by construction of the space $\rT$.  However, this is not true in general, and therefore it becomes important to estimate $E(\rT,K-1/2)$ for particular choices of reconstruction space $\rT$.
}

\subsection{Sampling at the critical density: the frame case} 
\label{s:frame_case} 
Unfortunately, the bound for $C(\Omega,\rT)$ declines as $\delta \rightarrow 1^{-}$, and is infinitely large at the critical value $\delta = 1$.  This result is sharp in the sense that there are countable nonuniform sampling schemes $\Omega = \{ \omega_n \}_{n \in \bbZ}$ (we now index over $\bbZ$ for convenience) with density $\delta = 1$ which are not complete (see \cite{ChristensenFramesAMS} or \cite{young}), and for which one therefore cannot expect stable or quasi-optimal reconstructions.  However, it is clear from considering uniform samples $\Omega = \{ n \}_{n \in \bbZ}$ that density $\delta = 1$ is permissible in some cases.  The standard approach to handle this ``critical'' density is to assume that the samples $\Omega = \{ \omega_n \}_{n \in \bbZ}$ give rise to an (unweighted) Fourier frame $\{ \E^{2 \pi \I \omega_n \cdot } \bbI_{[0,1]}(\cdot) \}_{n \in \bbZ}$ for $\rH$.  As we show next, stable reconstruction with NUGS is also possible in this setting.

\subsubsection{Background and notation}
Suppose that the ordered sequence $\{ \omega_n \}_{n \in \bbZ}$ gives rise to a Fourier frame for $\rH$.  In other words,
\bes{
A \| f \|^2 \leq \sum_{n \in \bbZ} | \hat{f}(\omega_n) |^2 \leq B \| f \|^2,\quad \forall f \in \rH,
}
for constants $0 < A \leq B < \infty$ (the frame constants).  
Note that the operator
\be{
\label{frame_operator}
\cS : \rH \rightarrow \rH,\ f \mapsto \sum_{n \in \bbZ} \hat{f}(\omega_n) \E^{2\pi \I \omega_n \cdot} \bbI_{[0,1]}(\cdot),
}
the so-called frame operator, is well-defined, linear, bounded and invertible, and satisfies
\bes{
A \| f \|^2 \leq \ip{\cS f}{f} \leq B \| f \|^2,\quad \forall f \in \rH.
}
Moreover, the truncated operators $\cS_N :  f \mapsto \sum^{N}_{n=-N} \hat{f}(\omega_n) \E^{2 \pi \I \omega_n \cdot} \bbI_{[0,1]}(\cdot)$ converge strongly to $\cS$ on $\rH$ as $N\rightarrow \infty$.

It shall be important later to have conditions under which a sequence $\{ \omega_n \}_{n \in \bbZ}$ gives rise to a Fourier frame.  Fortunately, in the one-dimensional setting, a near-characterization is known.  To state this, we first require several definitions:

\begin{enumerate}[(i)]
\item A sequence of points $\lambda_k \in \bbR$, $k \in I$, is called \textit{separated} if $| \lambda_k - \lambda_j | \geq \eta$, $j \neq k$, for some $\eta>0$. If $\{ \lambda_k \}_{k \in I}$ is a finite union of separated sets, then it is called a \textit{relatively} separated sequence.

\item For a sequence $\{ \omega_n \}_{n \in \bbZ}$, the lower Beurling density is defined by
\bes{
D^{-} = \lim_{r \rightarrow \infty} \frac{n^{-}(r)}{r},\qquad n^{-}(r) = \min_{t \in \bbR} \left | \left \{ n \in \bbZ : \omega_n \in (t,t+r) \right \} \right |.
}
\end{enumerate}
The following theorem, due to Jaffard \cite{Jaffard} and Seip \cite{SeipJFA}, gives an almost characterization of Fourier frames in terms of relative separation and the Beurling density:

\thm{
\label{t:frame_char}
If $\{ \omega_n \}_{n \in \bbN}$ is relatively separated and $D^{-} >1$ then $\{ \E^{2\pi \I \omega_n \cdot } \bbI_{[0,1]}(\cdot) \}_{n \in \bbZ}$ forms a frame for $\rH$.  Conversely, If $\{ \E^{2\pi \I \omega_n \cdot } \bbI_{[0,1]}(\cdot) \}_{n \in \bbZ}$ forms a frame for $\rH$ then $D^{-} \geq 1$ and $\{ \omega_n \}_{n \in \bbZ}$ is relatively separated.
}

Note that there exist both relatively separated sequences with $D^{-} =1$ which form frames and relatively separated sequences with $D^{-} =1$ which do not.  See \cite{ChristensenFramesAMS} for details.

\subsubsection{Stable reconstructions from frame samples} Let an ordered sequence $\{ \omega_n : n \in \bbZ \}$ give rise to a Fourier frame and let
$\Omega = \Omega_N = \{ \omega_n : |n|\leq N \}$.  According to Theorem \ref{t:bounds} stable reconstruction is possible with NUGS provided an admissible sampling operator exists.  Fortunately, this is always the case:

\thm{
\label{t:frame_admiss}
Let $\rT$ be a finite-dimensional subspace of $\rH$, and suppose that $\Omega_N = \{ \omega_n : |n|\leq N\}$, where $\{ \omega_n : n \in \bbZ \}$ gives rise to a Fourier frame.  Then the partial frame operator 
\be{
\label{partial_frame_op}
\cS_N : f \mapsto \sum^{N}_{n=-N} \hat{f}(\omega_n) \E^{2 \pi \I \omega_n \cdot},
}
is admissible for all sufficiently large $N$.  Specifically,
\be{
\label{frame_est}
C(\Omega,\rT) \leq  \frac{\sqrt{B}}{\sqrt{A - \tilde E(\rT,N)^2}},
}
where $A$ and $B$ are the frame constants and
\be{
\label{tilde_E}
\tilde E(\rT,N)^2 = \sup \left \{ \sum_{|n| > N} | \hat{f}(\omega_n) |^2 : \ f \in \rT, \| f \| =1 \right \}.
}
}
\prf{
The operator $\cS_N$ trivially satisfies conditions (i) and (ii) of Definition \ref{d:admissible}.  For the upper bound \R{approxParseval2} we merely note that $\ip{\cS_N f}{f} \leq \ip{\cS f}{f} \leq B \| f \|^2$, where $\cS$ is the frame operator \R{frame_operator}.  Moreover, since $\cS_N \rightarrow \cS$ strongly and $\rT$ is finite-dimensional, \R{approxParseval1} holds (with appropriate $C_1$) for all large $N$.  Specifically, for $f \in \rT$ we have 
\bes{
\ip{\cS_N f}{f} = \ip{\cS f}{f} - \ip{(\cS - \cS_N) f}{f} \geq A \| f \|^2 - \sum_{|n| > N} | \hat{f}(\omega_n) |^2 \geq \left (A - \tilde{E}(\rT,N)^2 \right ) \| f \|^2,
}
which gives $C_1(\Omega,\rT) \geq A - \tilde{E}(\rT,N)^2$.  We now apply the definition of $C$.
}

Note that this result is a trivial adaptation of results for GS proved in \cite{BAACHOptimality}.  We include it and its proof for completeness.  The novel results in the paper concerning frames come in the next two sections when we obtain estimates for the reconstruction constant $C(\Omega,\rT)$.

\subsection{Estimation of constants}\label{ss:estimation}
In summary, for reconstructions from nonuniform Fourier samples, we can now distinguish two cases.  When the samples $\Omega = \{ \omega_n : n=1,\ldots,N \}$ are $(K,\delta)$-dense, the results of \S \ref{ss:dense_stable} establish stable reconstruction with simple, numerically-verifiable, bounds for $C(\Omega,\rT)$.  Specifically, we may compute $C_1(\Omega,\rT)$ via an eigenvalue problem involving the matrix $A$ defined in \R{Abdef} (see \cite{BAACHOptimality} for details), and use the bound $C_2(\Omega) \leq (1 + \delta)^2$ obtained in Theorem \ref{t:Cbound_arbitrary} to give the computable estimate
$
C(\Omega,\rT) \leq C_{B}(\Omega,\rT) =  (1+\delta)/\sqrt{C_1(\Omega,\rT)}.
$
Conversely, if the samples are not $(K,\delta)$-dense, but arise from a Fourier frame, then as shown in the previous section, stable reconstruction is also possible.  Moreover, we have the estimate
$
C(\Omega,\rT) \leq C_{B}(\Omega,\rT) = \sqrt{B}/\sqrt{C_1(\Omega,\rT)},
$
in this setting.  Provided the upper frame bound $B$ is known, this estimate can be computed.  If $B$ is unknown (as is often the case in practice), then we may use a limiting process to compute $C_2(\Omega)$, and therefore $C(\Omega,\rT)$, to arbitrary accuracy. This is described in 
the following lemma:
\lem{\label{C_2_as_limit}
Suppose that $\Omega$ is finite and let $\cS:\rH \rightarrow \rH$ be a linear operator satisfying conditions $(\textnormal{i})$ and $(\textnormal{ii})$ of Definition $\ref{d:admissible}$.  Let $\rT_N$, $N \in \bbN$, be a sequence of finite-dimensional reconstruction spaces such that the corresponding orthogonal projections $\cP_N = \cP_{\rT_N}$ converge strongly to the identity on $\rH$.  Then
$
C_2(\Omega) = \lim_{N \rightarrow \infty} C_{3}(\Omega,\rT_N).
$
In particular, $C_2(\Omega)$ can be approximated to arbitrary accuracy by taking $N$ sufficiently large.
}
\prf{
Note first that $C_{3}(\Omega,\rT_N) \leq C_2(\Omega)$.  Let $f \in \rH$, $\| f \| =1$.  Then
\eas{
\ip{\cS f}{f} &= \ip{\cS \cP_N f}{\cP_N f} + \ip{\cS (f-\cP_N f)}{\cP_N f} +\ip{\cS f}{f- \cP_N f}
\\
& \leq C_{3}(\Omega,\rT_N) + 2 \sqrt{C_2(\Omega)} \sqrt{\ip{\cS( f - \cP_N f)}{f-\cP_N f}} .
}
Thus,
\eas{
C_3(\Omega,\rT_N) \leq C_2(\Omega) \leq C_3(\Omega,\rT_N) + 2 \sqrt{C_2(\Omega)} \sup_{f \in \rH, \| f \| = 1} \sqrt{\ip{\cS( f - \cP_N f)}{f-\cP_N f}}.
}
It suffices to show that the final term tends to zero as $N \rightarrow \infty$.  The operator $\cS$ is linear, bounded and, for any $g$, $\cS g$ depends only on the finite set of values $\{ \hat{g}(\omega) \}_{\omega \in \Omega}$.  Hence $\cS$ has finite rank. The result now follows from this and the strong convergence $\cP_N \rightarrow \cI$.
}

The key herein is that $C_3(\Omega,\rT_N)$ can easily be computed (see \cite{BAACHOptimality} for details).
Note that the same process can also be used in the case of $(K,\delta)$-dense samples.  But the improvement in doing so is likely marginal over the estimate $C_2(\Omega) \leq (1+\delta)^2 \leq 4$ (recall that $\delta \leq 1$).

\section{Reconstructions in wavelets}\label{s:wavelets}
We now consider the case of wavelet subspaces $\rT$.  Specifically, we answer the question of how large the bandwidth $K$ (or $N$ in the case of frame samples) needs to be to ensure stable, quasi-optimal reconstruction.  Our main result demonstrates $K$ (or $N$) need only scale linearly in $M = \dim(\rT)$ to guarantee this.  Thus, up to a constant, there is a one-to-one relationship between sampling bandwidth and wavelet scale.

\subsection{Preliminaries} \label{ss:wavelet_pre}
First, we recall the definition of a multiresolution analysis (MRA).
\begin{definition}
A multiresolution analysis of $\rL^2(\bbR)$ generated by a scaling function $\phi \in \rL^2(\bbR)$ is a nested sequence of closed subspaces $\{0\}\subseteq\cdots \subseteq V_{-2} \subseteq V_{-1} \subseteq V_{0} \subseteq V_1 \subseteq V_2 \subseteq \cdots \subseteq \rL^2(\bbR)$ such that
\begin{enumerate}
\item $\cup_{j \in \bbZ} V_j = \rL^2(\bbR)$ and $\cap_{j \in \bbZ} V_j = \{0\}$,
\item for all $j \in \bbZ$, $f(\cdot) \in V_j$ if and only if $f(2 \cdot) \in V_{j+1}$,
\item the collection $\{ \phi(\cdot - k) \}_{k \in \bbZ}$ forms a Riesz basis for $V_0$.
\end{enumerate}
\end{definition}
Also, recall that a system $\{ \phi(\cdot - k) \}_{k \in \bbZ}$ forms a Riesz basis for $V_0\subseteq\rL^2(\mathbb{R})$ if and only if there exists constants $d_1,d_2 > 0$ such that
\bes{
d_1 \sum_{k \in \bbZ} | \alpha_k |^2 \leq \nm{\sum_{k \in \bbZ} \alpha_k \phi(\cdot - k) }^2 \leq d_2 \sum_{k \in \bbZ} | \alpha_k |^2,\quad \forall \{ \alpha_k \}_{k \in \bbZ} \in l^2(\bbZ),
}
and $\{ \phi(\cdot - k) \}_{k \in \bbZ}$ forms an orthonormal basis for $V_0$ if and only if $d_1=d_2=1$.
This is equivalent to the condition
\be{
\label{RB_equivalent}
d_1 \leq \sum_{k \in \bbZ} \left | \hat{\phi}(k + \omega) \right |^2 \leq d_2,\quad a.e.\ \omega \in [0,1].
}
In particular, the optimal Riesz basis constants are given by
\bes{
d_1 =  \underset{\omega \in [0,1]}{\operatorname{essinf}} \sum_{k \in \bbZ}\left | \hat{\phi}(k + \omega) \right |^2 ,\qquad d_2 =  \underset{\omega \in [0,1]}{\operatorname{esssup}} \sum_{k \in \bbZ} \left | \hat{\phi}(k + \omega) \right |^2.
}
Now suppose that $\{ \psi_{j,k} \}_{j,k \in \bbZ}$ is a wavelet basis of $\rL^2(\bbR)$ associated to an MRA with a scaling function $\phi$. Our primary interest in this paper lies with wavelet bases with compact support. In this case, we can find a $p \in \bbN$ such that $\supp(\phi) \subseteq [-p+1,p]$.  Since we are interested in wavelet bases on the interval $[0,1]$,  following \cite{mallat09wavelet}, we consider three standard wavelet constructions -- periodic, folded and boundary wavelets -- with the corresponding MRA spaces 
\bes{
V^{\text{type}}_j = \spn \left \{ \phi^{\text{type}}_{j,k} :\ k=0,\ldots,2^{j}-1 \right \} ,\quad W^{\text{type}}_j = \spn \left \{ \psi^{\text{type}}_{j,k} : \ k=0,\ldots,2^{j}-1 \right \},
}
where $\text{type}\in\{\per,\fold,\ivl\}$ stands for periodic, folded or boundary wavelets, respectively. For a detailed construction of each basis type we refer to the appendix. 

Given $J \in \bbN_0$, we now introduce the finite-dimensional reconstruction space $\rT$ as follows
\be{
\label{wavelet_space_type}
\rT = V^{\text{type}}_J \oplus W^{\text{type}}_J \oplus W^{\text{type}}_{J+1} \oplus \cdots \oplus W^{\text{type}}_{R-1}.
}
Note that $\rT=V^{\text{type}}_R$ and also $\dim(\rT) = 2^R$. Since
\bes{
\supp (\phi_{R,k}) = [(k-p+1)/ 2^R , (k+p)/2^R ] \subseteq [0,1], \quad k=p,\ldots,2^R-p-1,
}
we have that $\phi^{\text{type}}_{R,k}(x) = \phi_{R,k}(x)$, $x \in [0,1]$, whenever $k=p,\ldots,2^R-p-1$.  Hence we may decompose the space $\rT$ into
\bes{
\rT = \rT^{\lt} \oplus \rT^i \oplus \rT^{\rt},
}
where $\rT^{i} = \spn \left \{ \phi_{R,k} : k=p,\ldots,2^R-p-1 \right \}$ contains unchanged interior  scaling functions with support in $(0,1)$ and 
\eas{
\rT^{\lt} &= \spn \left \{ \phi^{\text{type}}_{R,k} \bbI_{[0,1]} : k=0,\ldots,p-1 \right \}, \\
\rT^{\rt} &= \spn \left \{  \phi^{\text{type}}_{R,k}  \bbI_{[0,1]}  : k = 2^R-p,\ldots,2^R-1 \right \},
}
contain the modified scaling functions correspondingly for each basis type.  Here $\bbI_{[0,1]}$ is the indicator function of the interval $[0,1]$.  Whilst not strictly necessary at this point, we add this function to the definitions of $\rT^{\lt}$ and $\rT^{\rt}$ so as to clarify that they are to be considered as subspaces of $\rH = \{ g \in \rL^2(\bbR) : \supp(g) \subseteq [0,1] \}$ in our setting, and not $\rL^2(\bbR)$.

\rem{
\label{r:Why_not_DWT}
Periodic wavelet bases on $[0,1]$ are widely used in standard implementations of wavelets, since their construction is extremely simple.  However, the vanishing moments of the wavelet are lost due to the enforcement of periodic boundary conditions.  This effectively introduces a discontinuity of the signal at the boundaries, and translates into lower approximation orders $\textnormal{\cite{mallat09wavelet}}$. Folded wavelets remove the artificial signal discontinuity introduced by periodization and allow for one vanishing moment to be retained. This approach is most commonly used for the CDF wavelets $\textnormal{\cite{CDF}}$. However, since folded wavelets only retain one vanishing moment, they do not lead to high approximation orders for smooth functions.  To obtain such orders, one may follow the boundary wavelet construction, due to Cohen, Daubechies \& Vial $\textnormal{\cite{CDV}}$. These wavelets are particularly well suited for smooth functions.  Indeed, if $f \in \rH^{s}(0,1)$, where $\rH^{s}(0,1)$ denotes the usual Sobolev space and $0 \leq s < p$, then the error
\be{
\label{boundary_conv}
\| f - \cP_{\rT} f \| = \ord{2^{-s R}},\quad R \rightarrow \infty,
}
where $\rT$ is given by \R{wavelet_space_type} for $\textnormal{type}=\ivl$.  Since NUGS is quasi-optimal, we obtain exactly the same approximation rates when reconstructing $f$ from nonuniform Fourier samples, provided the bandwidth $K$ $($or $N$ in the frame case$)$ is chosen suitably large.  Corollary $\ref{c:super_NUGS}$ below establishes that $K$ $($or $N$$)$ need only scale linearly in $M = 2^R$ to guarantee this.  
}

\subsection{Main results}

Here, we state our main results regarding reconstructions in wavelet bases.  The proofs of these results are postponed to \S \ref{ss:proofs}.

\subsubsection{General wavelets}

We commence with the $(K,\delta)$-dense case:

\thm{
\label{t:main_Kdelta}
Let $\Omega$ be a $(K,\delta)$-dense sampling scheme and suppose that $\rT$ is the reconstruction space $\R{wavelet_space_type}$ of dimension $2^R$ generated by the first $2^R$ elements of a periodic wavelet basis $(\textnormal{type}=\per)$.  Suppose that either of the following conditions holds:
\begin{enumerate}[(i)]
\item the scaling function $\phi \in \rL^2(\bbR)$ and $\{ \phi(\cdot - k ) \}_{k \in \bbZ}$ forms an orthonormal basis of $V_0$,
\item the scaling function $\phi$ satisfies
\be{
\label{wavelet_smoothness}
| \hat{\phi}(\omega) | \leq \frac{c}{(1+|\omega|)^{\alpha}},\quad \omega \in \bbR,
}
for some $\alpha > \frac12$, and the system $\{ \phi(\cdot - k ) \}_{k \in \bbZ}$ forms a Riesz basis of $V_0$.
\end{enumerate}
Then for any $0 < \epsilon < 1- \delta$ there exists a $c_0 = c_0(\epsilon)$ such that if $K \geq c_0(\epsilon) 2^R$ then the reconstruction constant
\bes{
C(\Omega,\rT) \leq \frac{1+\delta}{1  - \delta- \epsilon}.
}
}

\thm{
\label{t:main_Kdelta_2}
Let $\Omega$ be a $(K,\delta)$-dense sampling scheme and suppose that either: $(\textnormal{i})$ $\rT$ is the reconstruction space $\R{wavelet_space_type}$ of dimension $2^R$ generated by the first $2^R$ elements of the folded wavelets basis $(\textnormal{type}=\fold)$, or $(\textnormal{ii})$ $\rT$ is the reconstruction space $\R{wavelet_space_type}$ of dimension $2^R$ generated by the first $2^R$ elements of the folded wavelets basis $(\textnormal{type}=\ivl)$.  Suppose that $\{ \phi(\cdot - k ) \}_{k \in \bbZ}$ is a Riesz basis for $V_0$ and that $\phi$ satisfies $\R{wavelet_smoothness}$ for some $\alpha > \frac12$.  Then given $0 < \epsilon < 1-\delta$ there exists a $c_0 = c_0(\epsilon)$ such that
\bes{
C(\Omega,\rT) \leq \frac{1+\delta}{1-\delta-\epsilon}, \quad  K\geq c_0(\epsilon) 2^R.
}
}

These theorems give the main result: bandwidth $K$ needs to scale linearly with $M = \dim(\rT) = 2^R$ in the case of wavelets.  Note that the smoothness assumption \R{wavelet_smoothness} is extremely mild.  For example, it holds if $\phi \in \rH^{\alpha}(\bbR)$ for $\alpha > \frac12$, and consequently includes all practical cases of interest. We remark also that the stipulation of a Riesz basis in these theorems is not necessary since this is implied by the MRA property.  It is included merely for clarity. 

We now give a similar result for the frame case:

\thm{
\label{t:main_frame}
Let $\Omega = \{ \omega_n : |n| \leq N\}$, where $\{ \omega_n : n \in \bbZ\}$ is a nondecreasing sequence that gives rise to a Fourier frame with frame bounds $A$ and $B$.  Let $\rT$ be the reconstruction space $\R{wavelet_space_type}$ of dimension $2^R$ consisting of either periodic, folded or boundary wavelets $(\textnormal{type}\in\{\per,\fold,\ivl\})$, and suppose that $\phi$ satisfies $\R{wavelet_smoothness}$ for some $\alpha > \frac12$.  Then given $0 < \epsilon < A$ there exists a $c_0 = c_0(\epsilon)$ such that
\bes{
C(\Omega,\rT)  \leq \sqrt{\frac{B}{A - \epsilon}}, \quad N \geq c_0(\epsilon) 2^R.
}
}

As explained in Remark \ref{r:Why_not_DWT}, boundary wavelets are an important case of these theorems.  Due to  \R{boundary_conv}, these results imply the following property of NUGS: up to constant factors, it obtains optimal convergence rates in terms of the sampling bandwidth when reconstructing smooth functions with boundary wavelets.  Specifically,
\cor{
\label{c:super_NUGS}
Let $\rT$ be the reconstruction space $\R{wavelet_space_type}$ consisting of the boundary wavelets $(\textnormal{type}=\ivl)$.  If $f \in \rH^s(0,1)$, where $0 \leq s < p$, let $\tilde{f}$ denote the NUGS reconstruction based on a sampling scheme $\Omega$.  Then $\| f - \tilde f \| = \ord{K^{-s}}$ if $\Omega$ is as in Theorem $\ref{t:main_Kdelta_2}$ and $\| f - \tilde{f} \| = \ord{N^{-s}}$ when $\Omega$ is as in Theorem $\ref{t:main_frame}$.
}

\subsubsection{Explicit estimates for Haar wavelets and digital signal models}\label{ss:exact_estimate}
Theorems \ref{t:main_Kdelta}--\ref{t:main_frame} do not give explicit bounds for the constant $C(\Omega,\rT)$.  In general, getting explicit estimates is difficult, due primarily to the contributions of the boundary subspaces $\rT^{\lt}$ and $\rT^{\rt}$.  However, for the case of Haar wavelets, there are no such terms, and this means that explicit bounds are possible.

One motivation for studying the Haar wavelet case is that it corresponds to the situation of a digital model for the signal $f$.  Specifically, the reconstruction space for Haar wavelets
\bes{
\rT = \spn \left \{  \phi  \cup \{ \psi_{j,k} : k=0,\ldots,2^j -1,\ j=0,\ldots,R-1 \} \right \},
}
is a special case corresponding to $M = 2^R$ of reconstruction space
\be{
\label{pcwseconstspace}
\rU = \rU_M = \left \{ g \in \rL^2(0,1) : g |_{[m/M,(m+1)/M)} = \mbox{constant},\  m=0,\ldots,M-1 \right \},
}
consisting of piecewise constant functions (i.e.\ digital signals where $1/M$ is the pixel size).  Note that
\be{
\label{UM_def}
\rU_M = \spn \left \{ \sqrt{M} \phi(M \cdot - m) : m=0,\ldots,M-1 \right \},
}
is a subspace generated by shifts of the pixel indicator function $\phi(x) = \bbI_{[0,1]}(x)$.  This digital signal model is popular in imaging.  In particular, it is the basis of the widely-used the fast, iterative reconstruction technique for MRI \cite{FesslerFastIterativeMRI} (see Remark also \ref{r:it_NUGS_relation}).

Our next result gives an explicit upper bound for the reconstruction constant $C(\Omega,\rT)$ in this case, and demonstrates that $C(\Omega,\rT)$ is mild whenever $M$ is at most $2K$.  
\thm{\label{t:Haar_sharp}
Let $\Omega$ be a $(K,\delta)$-dense sampling scheme and let $\rT\subseteq \rU_M$, where $\rU_M$ is given by $\R{UM_def}$ for $\phi(x) = \bbI_{[0,1]}(x)$ and $M \leq 2K$.  Then the following hold:
\begin{enumerate}[(i)]
\item If $2K/M \in \bbN$ then
\bes{
C(\Omega,\rT) \leq \frac{\pi}{2} \left ( \frac{1+\delta}{1-\delta} \right ).
}
\item If $2K/M \notin \bbN$ and $M \geq 2$ then
\bes{
C(\Omega,\rT) \leq c_0 \left ( \frac{1+\delta}{1-\delta} \right ),\quad c_0 = \frac{1}{\snc \left ( \pi/2 + \pi\delta  / M \right ) }. 
}
In particular, $c_0 \sim \pi/2$ for $M \gg 1$.
\end{enumerate}
}

\rem{
\label{r:it_NUGS_relation}
As noted previously, the well-known iterative reconstruction technique $\textnormal{\cite{FesslerFastIterativeMRI}}$ is a specific instance of NUGS corresponding to the choice $\R{pcwseconstspace}$ for $\rT$, where the term `iterative' refers to the use of conjugate gradient iterations combined with NUFFTs to solve the least squares problem $($recall the discussion in $\S \ref{ss:computation}$$)$.  Thus, Theorem $\textnormal{\ref{t:Haar_sharp}}$ provides an explicit guarantee for stable, quasi-optimal reconstruction with this method.
}

\subsection{Proofs}
\label{ss:proofs}

We first require the following lemma:

\lem{
\label{l:E_separation}
Let $I \subseteq \bbN$ be a finite index set and suppose that $\{ \varphi_n : n \in I \} \subseteq \rH$ is a Reisz basis for its span $\rT = \spn \{ \varphi_n : n \in I \}$ with constants $d_1$ and $d_2$.  Let $I$ be partitioned into disjoint subsets $I_1,\ldots,I_r$, and write $\rT_{i} = \spn \{ \varphi_n : n \in I_i \}$.  Let $E(\rT,z)$ and $\tilde{E}(\rT,N)$ be given by $\R{z_residual}$ and $\R{tilde_E}$ respectively.  Then  
\bes{
E(\rT,z) \leq \sqrt{\frac{d_2 }{d_1} \sum^{r}_{i=1} E(\rT_i,z)^2},\qquad \tilde E(\rT,N) \leq \sqrt{\frac{d_2 }{d_1} \sum^{r}_{i=1} \tilde E(\rT_i,N)^2}
}
}
\prf{
Let $f = \sum_{n \in I} \alpha_n \varphi_n \in \rT \backslash \{0\}$ and write $f = \sum^{r}_{i=1} f_i$, where $f_i = \sum_{n \in I_i} \alpha_n \varphi_n$.
Note that
\eas{
\| \hat{f} \|^2_{\bbR \backslash (-z,z)} \leq \left ( \sum^{r}_{i=1} \| \widehat{f_i} \|_{\bbR \backslash (-z,z)} \right )^2 \leq \left ( \sum^{r}_{i=1} E(\rT_i,z) \| f_i \| \right )^2 \leq \sum^{r}_{i=1} E(\rT_i,z)^2 \sum^{r}_{i=1} \| f_i \|^2.
} 
Also, since $\{ \varphi_n \}_{n \in I}$ forms a Riesz basis, we have $\sum^{r}_{i=1} \| f_i \|^2 \leq d_2/d_1 \| f \|^2$, and  therefore
\bes{
\| \hat{f} \|^2_{\bbR \backslash (-z,z)}/\|f\|^2 \leq d_2/d_1 \sum^{r}_{i=1} E(\rT_i,z)^2.
}
Taking the supremum over $f$ now gives the result for $E(\rT,z)$.  For $\tilde{E}(\rT,N)$, note first that
$\sum_{|n|>N}|\hat{f_i}(\omega_n)|^2< \infty$, $i=1,\ldots,r$, since $\{\omega_n\}_{n\in\bbZ}$ gives rise to the Fourier frame and  $f_i\in\rL^2(0,1)$, $i=1,\ldots,r$. Therefore, we  get 
$
\sqrt{\sum_{|n|>N}|\hat{f}(\omega_n)|^2} \leq \sum_{i=1}^{r} \sqrt{\sum_{|n|>N}|\hat{f_i}(\omega_n)|^2},
$
by using Minkowski's inequality.
Thus,
\bes{
\sum_{|n|>N}|\hat{f}(\omega_n)|^2 \leq  \left( \sum_{i=1}^{r} \tilde{E}(\rT^i,N) \| f_i \| \right)^2 \leq \frac{d_2}{d_1} \| f \|^2 \sum^{r}_{i=1} \tilde{E}(\rT^i,N)^2,
}
as required.
}

Recall that all the wavelet reconstruction systems introduced in \S \ref{ss:wavelet_pre} can be decomposed into interior wavelets having support in $[0,1]$ and boundary wavelets that intersect the endpoints $x=0,1$.  This lemma allows us to estimate the residuals $E(\rT,z)$ and $\tilde{E}(\rT,N)$ by considering each subspace separately. The next two propositions address the interior wavelets:

\prop{
\label{prop:linear_dense}
Let $\phi \in \rL^2(\bbR)$ have compact support and suppose that $\{ \phi(\cdot - k ) \}_{k \in \bbZ}$ forms a Riesz basis for its span with constants $d_1$ and $d_2$.  Let $M \in \bbN$, $M_1,M_2 \in \bbZ$ and
\bes{
\rT = \spn \left \{ \sqrt{M}\phi(M \cdot -m) :\ m=M_1,\ldots,M_2 \right \},
}
and suppose that $M,M_1,M_2$ are such that $\rT \subseteq \rH$.  Then
the following hold:
\begin{enumerate}
\item Given $\epsilon > 0$ there exists a $c_0 = c_0(\epsilon)$ such that
$
E(\rT,z)^2 < 1 - \frac{d_1}{d_2} + \epsilon,
$
for $z \geq c_0 M.$
\item Suppose that $\phi$ satisfies $\R{wavelet_smoothness}$ for some $\alpha > \frac12$.  Then there exists a $c_0 = c_0(\epsilon)$ such that
$
E(\rT,z)^2 < \epsilon,
$
for $z \geq c_0 M.$
\end{enumerate}
}

\prf{
Let $f \in \rT$ and write $f(x) = \sqrt{M} \sum^{M_2}_{k=M_1} a_k \phi(M x - k)$.
Since $\{ \phi(\cdot - k) \}_{k \in \bbZ}$ is a Riesz basis, we find that
\be{
\label{step00}
d_1 \sum^{M_2}_{k=M_1} | a_k |^2\leq \| f \|^2 \leq d_2 \sum^{M_2}_{k=M_1} | a_k |^2.
}
Moreover, a simple calculation gives that
\be{
\label{step0}
\hat{f}(\omega) =  \frac{1}{\sqrt{M}} \hat{\phi}\left (\frac{\omega}{M} \right ) \Psi \left (\frac{\omega}{M} \right ),\quad \omega \in \bbR,
}
where $\Psi(x) = \sum^{M_2}_{k=M_1}a_k \E^{-2\pi \I k x}$ is a trigonometric polynomial with 
$\| \Psi \|^2 = \sum^{M_2}_{k=M_1} | a_k |^2$.
Thus, by using \R{step00} we get 
\be{
\label{step1}
d_1 \| \Psi \|^2 \leq \| f \|^2 \leq d_2 \| \Psi \|^2.
}
We now estimate $\| \hat{f} \|^2_{(-z,z)}$.  By \R{step0}, we have
\eas{
\| \hat{f} \|^2_{(-z,z)} = \frac{1}{M} \int_{|\omega| < z} | \hat{\phi}(\omega/M)  |^2  \left | \Psi(\omega/M) \right |^2 \D \omega = \int_{|t| < z/M}  | \hat{\phi}(t) |^2 \left | \Psi(t) \right |^2 \D t.
}
Suppose that $z \geq M$ and write $\lfloor z/M \rfloor= n_0 +1$, where $n_0 \in \bbN_0$. Then
\bes{
\| \hat{f} \|^2_{(-z,z)} \geq \int^{n_0+1}_{t=-n_0} | \hat{\phi}(t) |^2 \left | \Psi(t) \right |^2 \D t  = \sum_{|n| \leq n_0} \int^{1}_{0} | \hat{\phi}(t+n) |^2 \left | \Psi(t+n) \right |^2 \D t.
}
Since $\Psi$ is $1$-periodic, and since  \R{step1} holds, we get
\bes{
\| \hat{f} \|^2_{(-z,z)} \geq \left (  \min_{t \in [0,1]}  \sum_{|n| \leq n_0}  | \hat{\phi}(n + t )|^2 \right ) \int^{1}_{0} | \Psi(t) |^2 \D t  \geq \frac{1}{d_2} \left ( \min_{t \in [0,1]} \sum_{|n| \leq n_0}  | \hat{\phi}(n + t )|^2 \right ) \| f \|^2. 
}
By  \cite[Lem.\ 5.4]{AHPWavelet}, there exists an $n_0 \in \bbN$ sufficiently large such that the term in brackets is greater than $d_1 - \epsilon d_2$.  Thus we get
$
\| \hat{f} \|^2_{(-z,z)} \geq \left ( d_1/d_2 - \epsilon \right ) \| f \|^2.
$
We now use the definition of $E(\rT,z)^2$ to complete part 1.\ of the proof.

Our approach for part 2.\ is similar, where we estimate the tail $\nmu{\hat{f}}^2_{\bbR \backslash (-z,z)}$.  Repeating the steps of the above proof, we find that
$
\| \hat{f} \|^2_{\bbR \backslash (-z,z)} \leq {1}/{d_1} \left ( \sup_{t \in [0,1]}  \sum_{|n| \geq n_0}  | \hat{\phi}(n + t )|^2 \right ) \| f \|^2.
$
Using the smoothness assumption \R{wavelet_smoothness}, we find that 
$
\sup_{t \in [0,1]} \sum_{|n| \geq n_0}  | \hat{\phi}(n + t )|^2 \lesssim (n_0)^{1-2 \alpha}.
$
Hence, if $z \geq c_0(\epsilon) M$ for some $c_0$, then $\| \hat{f} \|^2_{\bbR \backslash (-z,z)}  \leq \epsilon \| f \|^2$,
from which the result follows.
}

To prove our next result, we require the following useful lemma:

\lem{
\label{l:another_Grochenig}
Let $x_0 \leq x_1 < x_2 < \ldots < x_N \leq x_{N+1}$ where $N \in \bbN \cup \{ \infty \}$, and suppose that $\delta = \max_{n=0,\ldots,N} \{ x_{n+1} - x_n \} < \infty$.  Let $f \in \rH^1(a,b)$, where $a = \frac12(x_1+x_0)$, $b=\frac12(x_{N+1} + x_N)$ and $\rH^1(a,b)$ denotes the standard Sobolev space of first order on the interval $(a,b)$.  If $\mu_n = \frac12(x_{n+1}-x_{n-1})$, $n=1,\ldots,N$, then the following inequalities hold:
\bes{
\left ( \| f \|_{[a,b]} - \delta \| f' \|_{[a,b]} /\pi \right )^2 \leq  \sum^{N}_{n=1} \mu_n | f(x_n) |^2 \leq \left ( \| f \|_{[a,b]} + \delta \| f' \|_{[a,b]}/\pi \right )^2.
}
}
\prf{
The proof of this lemma is similar to that of Lemma \ref{l:FT_decay}.   Let $z_n = \frac{1}{2} (x_n + x_{n-1})$ and define $\chi(x) = \sum^{N}_{n=1} f(x_n) \bbI_{[z_n,z_{n+1})}(x)$.  Note that $z_1 = a$, $z_{N+1} = b$ and that $\sum^{N}_{n=1} \mu_n | f(x_n) |^2 = \| \chi \|^2_{[a,b]}$.  We now have 
$
\| f - \chi \|^2_{[a,b]} = \sum^{N}_{n=1} \int^{z_{n+1}}_{z_n} | f(x)  - f(x_n) |^2 \D x,
$
and after an application of Wirtinger's inequality, we obtain $\| f - \chi \|^2_{[a,b]} \leq \frac{\delta^2}{\pi^2} \| f' \|^2_{[a,b]}$.  This gives the result.
}

Having addressed the case of $(K,\delta)$-dense samples, we now  consider frame samples.  Recalling the setup of \S \ref{s:frame_case}, let $\{ \omega_n : n \in \bbZ \}$ be a nondecreasing sequence giving rise to a Fourier frame.  Set $\Omega_N = \{ \omega_n : |n| \leq N \}$, and suppose that $\cS_N$ is given by \R{partial_frame_op}.  

\prop{
\label{prop:frame_linear}
Let $\{ \omega_n \}_{n \in \bbZ} \subseteq \bbR$ be a nondecreasing sequence of frequencies that rise to a Fourier frame for $\rH$, and suppose that $\phi$ and $\rT$ are as in Proposition $\ref{prop:linear_dense}$.  If $\phi$ satisfies $\R{wavelet_smoothness}$ for some $\alpha > \frac12$, then given $\epsilon >0$ there exists a $c_0 = c_0(\epsilon)$ such that 
$
\tilde{E}(\rT,N) < \epsilon,
$
for all $N \geq c_0 M.$
}

For the proof of this result, we refer to the appendix.  We are now ready to prove  Theorems \ref{t:main_Kdelta} and \ref{t:main_Kdelta_2}. 

\prf{[Proof of Theorems $\ref{t:main_Kdelta}$ and $\ref{t:main_Kdelta_2}$]
By Theorem \ref{t:Cbound_arbitrary}, it suffices to consider $E(\rT,z)$.  Recall that in all three cases -- periodic, folded or boundary wavelets -- the reconstruction space $\rT$ can be decomposed as $\rT = \rT^{\lt} \oplus \rT^i \oplus \rT^{\rt}$.  Lemma \ref{l:E_separation} now gives
\bes{
E(\rT,z)^2 \leq \frac{d_2}{d_1} \left ( E(\rT^{\lt},z)^2 + E(\rT^i,z)^2 + E(\rT^\rt,z)^2 \right ).
}
The subspace $\rT^i$ contains wavelets supported in $[0,1]$, an application of Proposition \ref{prop:linear_dense} gives $E(\rT^i,z)^2 < \epsilon$ in both case (i) and case (ii) of Theorem  \ref{t:main_Kdelta} (recall in case (i) that $\{ \phi(\cdot-k) \}_{k \in \bbZ}$ is an orthonormal basis, and therefore $d_1=d_2 = 1$), as well as in Theorem  \ref{t:main_Kdelta_2}.  Thus it remains to show in all cases that $E(\rT^\lt,z)$ and $E(\rT^{\rt},z)$ can be made arbitrarily small with $z \gtrsim 2^R$

Consider the subspace $\rT^{\lt}$ (the case of $\rT^{\rt}$ is identical). For all three wavelet constructions, we may write
$
\rT^{\lt} = \spn \left \{ \Phi_{R,k} \bbI_{[0,1]} : k=0,\ldots,p-1 \right \},
$
where $\Phi_{R,k}$ is either $\phi^{\per}_{R,k}$ (periodic), $\phi^{\fold}_{R,k}$ (folded) or $\phi^{\ivl}_{R,k}$ (boundary).  The functions $\Phi_{R,k}\bbI_{[0,1]}$ form a Riesz basis for $\rT^{\lt}$ with bounds $d_1$ and $d_2$.  Hence, if 
$$
f = \sum^{p-1}_{k=0} \alpha_k \Phi_{R,k} \bbI_{[0,1]} \in \rT^{\lt},
$$
then $d_1 \sum^{p-1}_{k=0} | \alpha_k |^2 \leq \nm{f}^2 \leq d_2 \sum^{p-1}_{k=0} | \alpha_k |^2$.  Now consider $\nmu{\hat{f}}_{\bbR \backslash(-z,z)}$.  By the Cauchy--Schwarz inequality and the above inequality,
\bes{
\nmu{\hat{f}}_{\bbR \backslash(-z,z)} \leq \sum^{p-1}_{k=0} | \alpha_k | \nmu{(\Phi_{R,k} \bbI_{[0,1]})^{\wedge}}_{\bbR \backslash(-z,z)} \leq \sqrt{p / d_1} \| f \| \max_{0\leq k\leq p-1} \left \{ \nmu{(\Phi_{R,k} \bbI_{[0,1]})^{\wedge}}_{\bbR \backslash(-z,z)} \right \},
}
Thus, to complete the proof, we only need to show that there exists a $c_0 = c_0(\epsilon)$ such that
\be{
\label{conj}
\nmu{(\Phi_{R,k} \bbI_{[0,1]})^{\wedge}}_{\bbR \backslash(-z,z)} < \epsilon,\qquad \forall k=0,\ldots,p-1,
}
whenever $z \geq c_0(\epsilon) 2^R$.

Assume now that $2^{R-1} > p$.  Then one can determine the following:
\begin{itemize}
\item[(a)] For periodic wavelets, $\Phi_{R,k}(x) = \phi_{R,k}(x) + \phi_{R,k}(x-1)$.  
\item[(b)] For folded wavelets, $\Phi_{R,k}(x) = \phi_{R,k}(x) + \phi_{R,k}(-x)$.
\item[(c)] For boundary wavelets, $\Phi_{R,k}(x)$ can be written as a finite linear combination of the functions $\phi_{R,k}(x)$, where $k=-p+1,\ldots,p-1$.
\end{itemize}
Note that (a) and (b) follow by first writing $\phi^{\per}_{R,k}$ and $\phi^{\fold}_{R,k}$ in terms of infinite sums using the periodization and folding operations  
\bes{
f(x) \mapsto f^{\mathrm{per}}(x) = \sum_{k \in \bbZ} f(x+k),\quad f(x) \mapsto f^{\mathrm{fold}}(x) = \sum_{k \in \bbZ} f(x-2k) + \sum_{k \in \bbZ} f(2k-x),
}
and then by using the fact that $\supp(\phi) \subseteq [-p+1,p]$.  Case (c) was shown in \cite{CDV}.  Since in all cases $\Phi_{R,k}$ can be written as a finite sum with a number of terms independent of $R$, it therefore suffices to show that
\be{
\label{conj2}
\nmu{(\phi_{R,k} \bbI_{[0,1]})^{\wedge}}_{\bbR \backslash(-z,z)},\nmu{(\phi_{R,k}(\cdot-1) \bbI_{[0,1]})^{\wedge}}_{\bbR \backslash(-z,z)},  \nmu{(\phi_{R,k}(-\cdot) \bbI_{[0,1]})^{\wedge}}_{\bbR \backslash(-z,z)}< \epsilon,
}
where $k=-p+1,\ldots,p+1$ for the first term and $k=0,\ldots,p-1$ for the second two terms, whenever $z \geq c_0(\epsilon) 2^R$.  Note that
\bes{
\left | \left ( \phi_{R,k}(\cdot + l) \bbI_{[0,1]}\right )^{\wedge}(\omega) \right | = 2^{-R/2} \left | \int^{2^R(l+1)-k}_{2^R l-k} \phi(y) \E^{-2 \pi \I \omega y / 2^R} \D y \right |.
}
Suppose that $l=0$.  Then the integration interval is $[-k,2^R - k]$.  Since $\supp(\phi) = [-p+1,p]$, we can replace this by $[-k,p]$ to give 
\bes{
\left | \left ( \phi_{R,k}(\cdot) \bbI_{[0,1]}\right )^{\wedge}(\omega) \right | = 2^{-R/2} \left | \widehat{\phi^{[-k,p]}}(\omega/2^R) \right |,\quad k=-p+1,\ldots,p-1,
}
where $\phi^{[a,b]}(x) = \phi(x) \bbI_{[a,b]}(x)$ for $a<b$.  Similarly, for $l=-1$ we have
\bes{
\left | \left ( \phi_{R,k}(\cdot-1) \bbI_{[0,1]}\right )^{\wedge}(\omega) \right | = 2^{-R/2} \left | \widehat{\phi^{[-p+1,k]}}(\omega/2^R) \right |,\quad k=0,\ldots,p-1.
}
Likewise
\bes{
\left | \left ( \phi_{R,k}(-\cdot) \bbI_{[0,1]}\right )^{\wedge}(\omega) \right | = 2^{-R/2} \left | \widehat{\phi^{[-p+1,k]}}(-\omega/2^R) \right |,\quad k=0,\ldots,p-1.
}
Thus, to establish \R{conj2}, and therefore \R{conj}, it suffices to estimate the Fourier transforms of the functions $\phi^{[a,b]}$ for $(a,b) = (-k,p)$, $k=-p+1,\ldots,p-1$, and $(a,b) = (-p+1,k)$, $k=0,\ldots,p-1.$  We now note the following:
$
\nmu{2^{-R/2} f(\cdot / 2^R ) }_{\bbR \backslash (-z,z)} = \nmu{f}_{\bbR \backslash ( -z/2^R , z/2^R)},\quad f \in \rL^2(\bbR).
$
In particular, for any fixed $f$, 
\be{
\label{birds}
\nmu{2^{-R/2} f(\cdot / 2^R ) }_{\bbR \backslash (-z,z)} < \epsilon,
}
provided $z \geq c 2^R$ for appropriately large $c > 0$.  Since the total number of functions $\phi^{[a,b]}$ is less than $2p$, and hence bounded independently of $R$, we obtain \R{conj2} and thus \R{conj}.
}

\prf{[Proof of Theorem $\ref{t:main_frame}$]
By Theorem \ref{t:frame_admiss}, we may consider $\tilde{E}(\rT,N)$.  Proceeding in a similar manner to the previous proof, we see from Lemma \ref{l:E_separation} that it suffices to estimate $\tilde{E}(\rT^{i},N)$, $\tilde{E}(\rT^{\lt},N)$ and $\tilde E(\rT^{\rt},N)$ separately.  As before, $\tilde{E}(\rT^{i},N)$ can be bounded using Proposition \ref{prop:frame_linear}, and hence it remains to derive bounds for $\tilde{E}(\rT^{\lt},N)$ and $\tilde E(\rT^{\rt},N)$ only. If we now argue in an identical way to the previous proof, i.e.\ by writing the spaces $\rT^{\lt}$ and $\rT^{\rt}$ as linear combinations of the functions $\phi^{[a,b]}$ whose total number is independent of $R$, then we see that it suffices to show the following: for an arbitrary function $f \in \rL^2(0,1)$, 
\be{
\label{conj3}
2^{-R} \sum_{|n| > N} \left | \hat{f}(\omega_n/2^R) \right |^2 < \epsilon,
}
provided $N \geq c 2^R$ for some $c>0$ depending only on $f$ (this replaces the condition \R{birds} in the previous proof).
Recall from the proof of Proposition \ref{prop:frame_linear} that we may assume without loss of generality that the frame sequence $\{ \omega_n \}_{n \in \bbZ}$ is separated with separation at least $\eta/2$ and maximal spacing at most $\eta$.  Thus the points $\{\tilde{\omega}_n \}_{n \in \bbZ}$, where $\tilde{\omega}_n = \omega_n/2^R$, have maximal spacing at most $\eta / 2^R$ and we find that
$
2^{-R} \sum_{|n| > N} \left | \hat{f}(\omega_n/2^R) \right |^2 \leq \frac{2}{\eta} \sum_{|n| > N} \mu_n | \hat{f}(\tilde \omega_n) |^2,
$
where $\mu_n = \frac{\tilde{\omega}_{n+1} - \tilde{\omega}_{n-1}}{2}$.  Since $f \in \rH$ we may apply Lemma \ref{l:another_Grochenig} to get
\bes{
2^{-R} \sum_{|n| > N} \left | \hat{f}(\omega_n/2^R) \right |^2 \leq \frac{2}{\eta} \left [  \left ( \|\hat{f}\|_{J_{+}} + \frac{\eta}{2^R \pi} \|\hat{f}'\|_{J_{+}} \right )^2 + \left ( \|\hat{f}\|_{J_{-}} + \frac{\eta}{2^R \pi} \|\hat{f}'\|_{J_{-}} \right )^2 \right ],
}
where $J_+ = (\tilde \omega_N , \infty)$ and $J_{-} = (-\infty,\tilde \omega_{-N} )$.  To obtain \R{conj3} we merely note that $\hat{f}' = \widehat{f_1}  \in \rL^2(\bbR)$, where $f_1(x) = x f(x)$, and $\max \{ \tilde{\omega}_N, -\tilde{\omega}_{-N} \}\gtrsim N / 2^R$ for large $N$.
}

Finally, we prove Theorem \ref{t:Haar_sharp}, which gives an explicit upper bound for the reconstruction constant in the case of reconstructing in Haar wavelets:

\prf{[Proof of Theorem $\ref{t:Haar_sharp}$]
Since we have already shown have $C_2(\Omega)\leq (1+\delta)^2$, and since $\rT\subseteq \rU_M$, it is enough to estimate $C_1(\Omega,\rU_M)$. For any $f\in\rU_M$, we can write
$
f(x)=\sqrt{M}\sum_{m=0}^{M-1}a_m\phi(Mx-m).
$
Therefore, as before, we get
$\hat{f}(\omega)=\frac{1}{\sqrt{M}} \hat{\phi}\left(\frac{\omega}{M}\right) \tilde{\Psi}\left(\frac{\omega}{M}\right)$,
where, for $M_0 = \lfloor M/2 \rfloor$,
\bes{
\tilde{\Psi}(x) = \sum^{M-1}_{m=0} a_m \E^{-2 \pi \I m x} = \E^{-2 \pi \I M_0 x} \sum^{M-M_0-1}_{m=-M_0} a_{m+M_0} \E^{-2 \pi \I m x} = \E^{-2 \pi \I M_0 x} \Psi(x),
}
and $\Psi(x) =  \sum^{M-M_0-1}_{m=-M_0} a_{m+M_0} \E^{-2 \pi \I m x}$.  Note that $\Psi$ is a trigonometric polynomial of degree at most $M_0$ and moreover, since $\{ \phi(\cdot - k ) \}_{k \in \bbZ}$ is an orthonormal basis, we have $\| \Psi \|^2 = \| f \|^2$.  Set $x_n = \omega_n / M$ for $n=0,\ldots,N+1$ and let $\nu_n = \frac12(x_{n+1} - x_{n-1})$.  Then we have
\be{
\label{orange}
\ip{\cS f}{f} = \sum^{N}_{n=1} \nu_n | \Psi(x_n) |^2 | \hat{\phi}(x_n) |^2 .
}
Let us first consider the case $2K/M \in \bbN$.  Note that $\rU_M \subseteq \rU_{2K}$ in this case, and therefore it suffices to prove the result for $M=2K$.  After an application of Lemma \ref{l:another_Grochenig}, we obtain
\bes{
\ip{\cS f}{f} \geq \min_{n=1,\ldots,N} | \hat{\phi}(x_n) |^2  \left ( \| \Psi \|_{[a,b]} - \frac{\delta}{2 K \pi} \| \Psi' \|_{[a,b]} \right )^2 \geq d_0 \left ( \| \Psi \|_{[a,b]} - \frac{\delta}{2 K \pi} \| \Psi' \|_{[a,b]} \right )^2,
}
where $a = \frac12(x_1 + x_0) = \frac12(x_1+x_N) - \frac12$, $b = \frac12(x_N + x_{N+1}) = \frac12 (x_1 + x_N) + \frac12$ and $d_0 = \min_{\omega \in [-1/2,1/2]} | \hat{\phi}(\omega) |^2$.  Note that the second inequality here follows from the observation that $|x_n| = |\omega_n | / M \leq K/M \leq 1/2$ since the frequencies $\omega_n$ are $(K,\delta)$-dense.
Since $b-a = 1$ and $\Psi$ is periodic, we therefore have
\bes{
\ip{\cS f}{f} \geq d_0 \left ( \| \Psi \| - \frac{\delta}{2 K \pi} \| \Psi' \| \right )^2 \geq d_0 \left ( 1 - \frac{\delta M_0}{K} \right )^2 \| \Psi \|^2 \geq d_0 \left ( 1 - \delta \right )^2 \| \Psi \|^2,
}
where the penultimate inequality follows from $\| \Psi' \| \leq 2 M_0 \pi \| \Psi \|$.  To complete the proof, we note that $| \hat{\phi}(\omega) | = | \snc(\omega \pi) |$ and that $|\snc(\omega \pi)| \geq |\snc(\pi/2)| = 2/\pi$ for $\omega \in [-1/2,1/2]$.

Now suppose that $M \leq 2K$ is arbitrary.   In this case, our first step is to introduce a new subset of points $\{ \tilde{x}_p \}^{\tilde{N}}_{p=1}$.  We do this as follows.  Let $n'$ be the largest $n$ such that $x_n \leq -1/2$, and let $n''$ be the smallest $n$ such that $x_n \geq 1/2$.  If $\tilde{N} = n'' - n'+1$, let
$
\tilde{x}_p = x_{p+n'-1},
$
$p=0,\ldots,\tilde{N},
$
and 
$\tilde{x}_{\tilde{N}+1} = 2+x_{n'-1} +x_{n'}- x_{n''}.$
Let $\tilde{\nu}_p = \frac12 ( \tilde{x}_{p+1} - \tilde{x}_{p-1})$, and note that
$\tilde{\nu}_{p} = \nu_{p+n'-1}$ for $p=1,\ldots,\tilde{N}-1$.
Moreover, by definition of $n'$ and $n''$, we have
\bes{
\tilde{\nu}_{\tilde{N}} = \frac12 \left (\tilde{x}_{\tilde{N}+1} - \tilde{x}_{\tilde{N}-1} \right ) = \frac12 \left ( 2 + x_{n'-1} + x_{n'} - x_{n''} - x_{n''+1} \right ) + \nu_{n''} \leq \nu_{n''}.
}
Therefore, we now obtain the following from \R{orange}:
\bes{
\ip{\cS f}{f} \geq \min_{n=n',\ldots,n''} | \hat{\phi}(x_n) |^2 \sum^{\tilde{N}}_{p=1} \tilde{\nu}_p | \Psi(\tilde{x}_p) |^2.
}
Since the frequencies $\omega_n$ are $(K,\delta)$-dense, we have that $x_{n'} \geq - 1/2 - \delta/M$ and $x_{n''} \leq 1/2 + \delta/M$.  This and an application of Lemma \ref{l:another_Grochenig} now give
\bes{
\ip{\cS f}{f} \geq d_0 \left ( \| \Psi \|_{[a,b]} - \frac{\delta}{M \pi} \| \Psi'\|_{[a,b]} \right )^2,\qquad d_0 = \min_{\omega \in[-1/2-\delta/M,1/2+\delta/M]} | \hat{\phi}(\omega) |^2,
}
where
$
a = \tfrac12 \left ( \tilde{x}_1 + \tilde{x}_0 \right ) = \tfrac12 \left (x_{n'} + x_{n'-1} \right ),
$
and
$
b = \tfrac12 \left ( \tilde{x}_{\tilde{N}} + \tilde{x}_{\tilde{N}+1} \right )  = a + 1.
$
Since $|b-a| = 1$, we now argue exactly as before to give the first result.
}

\section{Bandwidth and ill-conditioning}\label{s:bandwidth}

Having shown that stable reconstruction is possible provided the bandwidth $K$ scales linearly with the dimension $M = 2^R$ of the wavelet reconstruction space, we now consider the constant of this scaling:

\thm{
\label{t:exp_blowup}
Let $\Omega = \{ \omega_n : n =1,\ldots,N \} \subseteq [-K,K]$ for some $K > \frac{2}{\pi^2} + \frac12$ and suppose that $\cS$ is given by $\R{weightedFS}$ with weights $\R{weights}$.  Let $\rT$ be the reconstruction space corresponding to either periodic, folded  or boundary wavelets, where $2^{R-1} > K$.  Then the reconstruction constant satisfies
\bes{
C(\Omega,\rT) \geq c_1 / \sqrt{K} \exp \left ( c_2 (1-z) 2^R \right ),
}
where $z = \max \{ \frac12 , K/2^{R-1} \}$ and $c_1,c_2>0$ depend only on $\phi$.
}

This theorem, which generalizes a result proved in \cite{AHPWavelet} to the case of nonuniform samples, establishes the following.  Suppose that the size $M = 2^R$ of the reconstruction space is roughly $2\alpha K$.  If $\alpha > 1$ then the reconstruction constant $C(\Omega,\rT)$ blows up exponentially fast as $M \rightarrow \infty$.  In other words, if the bandwidth $K$ of the sampling is not sufficiently large in comparison to the wavelet scale $R$, then ill-conditioning is necessarily witnessed in the reconstruction.  Note that this theorem does not assume density of the samples, just that their maximal bandwidth is $K$.  In particular, even if $\hat{f}(\omega)$ were known for arbitrary $| \omega | \leq K$ one would still have the same result, i.e.\ insufficient sampling bandwidth implies ill-conditioning.  

It is instructive to compare this result with Theorem \ref{t:Haar_sharp}, which estimates the reconstruction constant for Haar wavelets.  If $M \approx 2 \alpha K$ then Theorem \ref{t:Haar_sharp} demonstrates that $C(\Omega,\rT)$ is bounded whenever $\alpha$ is less than or equal to the critical value $\alpha_0 = 1$.  Conversely, if $\alpha > \alpha_0$ then exponential ill-conditioning necessarily results as a consequence of Theorem \ref{t:exp_blowup}.  For other wavelets, Theorems \ref{t:main_Kdelta} and \ref{t:main_Kdelta_2} show that stable reconstruction is possible for sufficiently small scaling $\alpha$, but unlike the Haar wavelet case, they do not establish the exact value for $\alpha_0$ that delineates the stability and instability regions.  

Theorem \ref{t:exp_blowup} follows immediately from the following lemma:
\lem{
\label{l:exp_blowup}
Let $\Omega$ and $\cS$ be as in Theorem $\ref{t:exp_blowup}$.  Let $\rT \subseteq \rH$ and suppose that $\rT \supseteq \rU$, where
\bes{
\rU = \spn \left \{ \sqrt{M} \phi(M \cdot - m) : m=M_1,\ldots,M_2 \right \}.
}
for some $M \in \bbN$, $M_1,M_2 \in \bbZ$ and $M > 2K$.  If $\{ \phi(\cdot-k) \}_{k \in \bbZ}$ is a Riesz basis for its span with bounds $d_1$ and $d_2$ then
\bes{
C(\Omega,\rT) \geq c_1  \sqrt{\frac{d_1}{d_2 K}} \exp \left [ c_2 (M_2 - M_1 -2 ) (1-z) \right ],
}
where $z = \max \{ \frac12 , 2K/M \}$, and $c_1,c_2>0$ are independent of $\Omega,K,M,M_1,M_2$ and $\phi$.
}

\prf{[Proof of Theorem \ref{t:exp_blowup}]
In each case, we merely set $\rU = \rT^i$ to be the space spanned by the interior wavelets.  The result follows immediately from Lemma \ref{l:exp_blowup}.
}

To prove Lemma \ref{l:exp_blowup}, we require the following result (see \cite[Prop.\ 6.2]{AHPWavelet} for a proof):
\lem{
\label{l:trig_exp_blowup}
Let $P \in \bbN$ and $z \in (0,1/2)$.  Then there exists a constant $c>0$ independent of $P$ and $z$ such that, if $z' = \max \{ 1/4 , z \}$, then
\bes{
\sup \left \{ \frac{\sup_{|t| \leq 1/2} | \Psi(t) | }{\sup_{|t| \leq z } | \Psi(t) | } : \Psi(t) = \sum_{|n| \leq P} a_k \E^{\I 2\pi k t},\  a_k \in \bbC \right \} \geq \exp \left( c P (1/2-z') \right).
}
}

\prf{[Proof of Lemma $\ref{l:exp_blowup}$]
Note that $C(\Omega,\rT) \geq C(\Omega,\rU)$.  Let $f \in \rU$.  Then
\bes{
\ip{\cS f}{f} = \frac{1}{M} \sum^{N}_{n=1} \mu_n | \hat{\phi}(\omega_n/M) |^2 | \Psi(\omega_n/M) |^2,
}
where  $\Psi(x) = \sum^{M_2}_{k=M_1}a_k \E^{-2\pi \I k x}$  satisfies $d_1 \| \Psi \|^2 \leq \| f \|^2 \leq d_2 \| \Psi \|^2$.  Thus
\eas{
\ip{\cS f}{f} &\leq \sup_{| \omega | \leq K/M } | \hat{\phi}(\omega) |^2 \sup_{|t| \leq K/M} | \Psi(t) |^2 \left ( \frac{1}{M} \sum^{N}_{n=1} \frac{\omega_{n+1}-\omega_{n-1}}{2} \right ) 
\\
&= \frac{2 K}{M} \sup_{| \omega | \leq K/M } | \hat{\phi}(\omega) |^2 \sup_{|t| \leq K/M} | \Psi(t) |^2 \leq \frac{2 K d_2}{M} \sup_{|t| \leq K/M} | \Psi(t) |^2,
}
where the final inequality follows from \R{RB_equivalent}.  The definition \R{approxParseval1} of $C_1(\Omega,\rU)$, now gives
\bes{
C_1(\Omega,\rT) \leq C_1(\Omega,\rU) \leq \frac{2 K d_2}{M d_1}  \inf_{\Psi \in \rV} \left \{ \frac{\sup_{|t| \leq K/M} | \Psi(t) |^2}{\| \Psi \|^2 }\right \},
}
where
$
\rV = \left \{ \sum^{M_2-M_3}_{k=M_1-M_3} a_k \E^{ 2 \pi \I k x} : a_k \in \bbC \right \},\ M_3 = \left \lceil \frac{M_1+M_2}{2} \right \rceil.
$
Since $M_2-M_1 \leq M$ we have $| \Psi(t) |^2 \leq (M+1) \| \Psi \|^2$, and therefore
\be{
\label{exp_step1}
C_1(\Omega,\rT) \leq \frac{d_2}{d_1} (2K+1) \inf_{\Psi \in \rV} \left \{ \frac{\sup_{|t| \leq K/M} | \Psi(t) |^2}{\sup_{|t| \leq 1/2} | \Psi(t) |^2 }\right \}.
}
We shall return to this in a moment.  First, let us consider $C_2(\Omega)$.  We wish to show that $C_2(\Omega) \geq c$ for any $\Omega$ for some $c>0$.  Suppose that $\Omega$ is $(K,\delta)$-dense.  Then by Lemma \ref{l:FT_decay},
$
C_2(\Omega) \geq  \left ( \sqrt{1- \nmu{\hat{f}}^2_{\bbR \backslash I} /  \nm{f}^2  } - \delta \right )^2,
$
where $I = (-K+\frac12 \delta,K-\frac12 \delta)$.  Let $f(x) = \bbI_{[0,1]}(x)$, so that $\hat{f}(\omega) = \E^{-\I \pi \omega} \snc (\omega \pi )$.  Then
$
C_2(\Omega) \geq  \left ( \sqrt{1-2/(\pi^2 ( K-1/2))} - \delta \right )^2.
$
Now suppose that $\Omega$ is not $(K,\delta)$ dense.  Then there exists an $n=1,\ldots,N$ such that $\omega_{n+1} - \omega_n \geq \delta$, and therefore $\mu_n \geq \delta/2$.  Hence
\bes{
C_2(\Omega) \|f\|^2 \geq \ip{\cS f}{f} \geq \tfrac12 \delta | \hat{f}(\omega_n) |^2,\quad \forall f \in \rH.
}
Picking $f(x) = \E^{2\pi\I  \omega_n  x}$, we therefore obtain $C_2(\Omega) \geq \frac12\delta $.  Hence in general
\bes{
C_2(\Omega) \geq \max \left \{ \delta/2 ,  \left ( \sqrt{1-2/(\pi^2(K-1/2))} - \delta \right )^2 \right \},\quad \forall \delta \in (0,1).
}
Since $1-2/(\pi^2 ( K-\frac12)) > 0$, and since $\delta > 0$ was arbitrary, we now find that $C_2(\Omega) \geq c^2$ for any $\Omega$.
Combining this with \R{exp_step1}, we now find that
\bes{
C(\Omega,\rT) \geq c \sqrt{\frac{d_1}{d_2 K}}\sup_{\Psi \in \rV} \left \{ \frac{\sup_{|t| \leq 1/2} | \Psi(t) |}{\sup_{|t| \leq K/M} | \Psi(t) |} \right \}.
}
To complete the proof, we first note that
$
\min \{ M_2 - M_3 , M_3 - M_1 \} \geq (M_2-M_1-1)/2.
$
Thus, $\rV$ contains all trigonometric polynomials of degree $ \left\lfloor\frac{M_2-M_1-1}{2}\right\rfloor \geq \frac{M_2-M_1}{2} -1$.  An application of Lemma \ref{l:trig_exp_blowup} now gives the result.
}

\rem{
\label{r:optimality}
Using techniques of $\textnormal{\cite{BAACHOptimality}}$, one can show that the reconstruction constant $C(\Omega,\rT)$ is essentially universal.  Specifically, any reconstruction algorithm that is so-called \textit{perfect} must have a condition number that is at least that of NUGS.  In particular, noting Theorem $\ref{t:exp_blowup}$, we see that to recover wavelet coefficients up to scale $R$ stably and accurately, it is necessary to take samples from a bandwidth $K$ that is at least $2^{R-1}$, regardless of the method used. 
}

\section{Numerical examples}\label{s:num_exp}

We now present several numerical examples to illustrate the NUGS.  We will focus on the following three nonuniform sampling schemes:
\begin{enumerate}[(i)]
 \item{Jittered sampling:} Let $K>0$ and $\eta,\epsilon \in (0,1)$ be such that $\epsilon + 2 \eta < 1$. Set $\tilde{N} = \left \lfloor \frac{K}{\epsilon} \right \rfloor$ and $N=2\tilde{N}+1$. The jittered sampling scheme is given by
$
\Omega_{N} = \{ \omega_{1},\ldots,\omega_{N} \}, 
$
where $\omega_n=n\epsilon + \eta_n,$ $n=-\tilde{N},\ldots,\tilde{N},$
and $\eta_n \in (-\eta,\eta)$ is chosen uniformly at random.  Note that $\Omega_{N}$ is $(K+\eta,\delta)$-dense, where $\delta=\epsilon+2 \eta$.  This sampling scheme is a standard model for jitter error in MRI caused by the measurement device not scanning exactly on a uniform grid \cite{GrochenigStrohmerMarvasti}. 

 \item{Log sampling:} Let $K>0$, and let $\nu$ and $\delta$ be fixed parameters such that $2 \times 10^{-\nu} < \delta$. Set $\tilde{N} = \left \lceil - \frac{\log_{10} K + \nu}{\log_{10}(1-\delta/K) } \right \rceil$ and  $N=2(\tilde{N}+1)$. Log sampling scheme is given by
\bes{
\Omega_N = \{ - \omega_n \}^{\tilde{N}}_{n=0} \cup \{ \omega_n \}^{\tilde{N}}_{n=0},\quad\text{where}\quad \omega_n = 10^{-\nu + \frac{n}{\tilde{N}} (\log_{10} K + \nu ) },\quad n=0,\ldots,\tilde{N}.
}
Note that this gives a $(K,\delta)$-dense sampling sequence.  This sampling scheme is a one-dimensional model for a two- or three-dimensional spiral sampling trajectory.  Such trajectories are popular in MRI applications (see \S \ref{s:introduction}).
  \item{Seip's frame:}  For a given $N\in\bbN$, set
$
\Omega_N = \{ \omega_n \}^{-N}_{n=-1} \cup \{ \omega_n \}^{N}_{n=1},
$
where $\omega_n = n (1-|n|^{-1/2})$ and $|n| \geq 1.$
In \cite{SeipJFA}, it is shown that the infinite set of frequencies $\Omega = \Omega_{\infty}$ gives rise to a Fourier frame with density $\delta=1$.
\end{enumerate}

The main result proved in \S \ref{s:wavelets} is that one requires a linear scaling of the bandwidth $K$ or truncation index $N$ with the parameter $M = 2^R$ for stable reconstruction in wavelet subspaces.  This is illustrated in Table \ref{tab:linear_scale} for the Haar and DB4 wavelets.  Note that the constant of the scaling is roughly $1/2$, i.e.\ $K$ (or $N$) behaves like $\beta 2^R$ with $\beta\approx 1/2$.  In the case of Haar wavelets, this is due to the explicit estimates of Theorem \ref{t:Haar_sharp}.

\begin{table}[H]
\centering
\scalebox{0.84}{
\begin{tabular}{ |c|c|c|cccccc||c|c|c|cccccc| }
\hline
$\rT$                                &$\Omega$&$2^R$&  32 & 64 & 128 & 256 &  512 & 1024 & $\rT$ & $\Omega$ &  $2^R$ &  32 & 64 & 128 & 256 &  512 &1024 \\ \hline
\multirow{2}{*}{\rotatebox{0}{Haar}} & Log   & $K$ &  16 & 32 & 64 & 128 &  256 & 512 & \multirow{2}{*}{\rotatebox{0}{DB4}} & Log & $K$ &  16 & 32 & 64 & 128 &  256 & 512  \\ 
                                     & Frame &  $N$  & 20 & 38 & 72 & 139 &  272 & 535 & & Frame &  $N$  & 20 & 38 & 72 & 139 &  272 & 535 \\ \hline
\end{tabular}
}
\caption{\small{For a given number of reconstruction vectors $2^R$, the smallest value of $K$ (or $N$) is shown such that the reconstruction constant $C(\Omega,\rT)$ is at most $100$, where the reconstruction constant is estimated by using the results given in \S \ref{ss:estimation}. This is done for different reconstruction spaces $\rT$ -- Haar and DB4 -- and  for different sampling schemes $\Omega$: Seip's frame sequence and log sampling scheme with $\delta=0.95$ and $\nu=0.33$.} 
}
\label{tab:linear_scale}
\end{table}

\begin{table}[H]
\centering
\scalebox{0.85}{
\begin{tabular}{ |c|c|cccccc| }
\hline
\multirow{2}{*}{\rotatebox{0}{$\rT$}} &  $c_0$ &  0.3125 & 0.3750 & 0.4375 & 0.5000 &  0.5625 &  0.6250  \\
& $K$ &  20 & 24 & 28 & 32 & 36 & 40    \\  \hline
\multirow{2}{*}{\rotatebox{0}{Haar}} & $\kappa(A)$ &  5.8569e15 & 2.9255e12 &	1.8347e05 & 1.7835 &  1.6474  & 1.5768  \\ 
&  $\frac{\| f - \tilde{f} \|}{\| f - \cP_\rT f \|}$  & 8.6294e04 &  7.3412e04 &  14.4886 &  1.0016 &  1.0016 &  1.0016 \\ \hline
\multirow{2}{*}{\rotatebox{0}{DB4}} & $\kappa(A)$ &  5.0079e15 & 2.6583e12 & 1.2918e05 & 1.6126 & 1.4744 & 1.4355  \\ 
&  $\frac{\| f - \tilde{f} \|}{\| f - \cP_\rT f \|}$  & 4.0459e06	 & 3.2764e06 & 303.3421 & 1.0013 & 1.0009 & 1.0008  \\ \hline
\end{tabular}
}
\caption{\small{The condition number $\kappa(A)$ and the error ${\| f - \tilde{f} \|}/{\| f - \cP_\rT f \|}$ are shown for different bandwidths $K=c_02^R$ and different reconstruction spaces: Haar and DB4 wavelets, where $2^R=64$ is taken. The jittered sampling scheme is used for $\epsilon=0.6$ and $\eta=0.15$, and the function $f(x)=1/2 \cos(4 \pi x)$ is tested.
}}
\label{tab:blowup}
\end{table}

Theorem \ref{t:exp_blowup} provides a lower estimate for such scaling.  In particular, if the scaling $\beta$ is less than $1/2$ then exponential instability necessarily results in the reconstruction, regardless of the wavelet basis used.  This is shown in Table \ref{tab:blowup} for both Haar and DB4 wavelets.  Note also that in the unstable regime, i.e.\ $\beta < 1/2$, the reconstruction $\tilde{f}$ is also far from quasi-optimal.

\begin{table}[H]
\centering
\scalebox{0.84}{
\begin{tabular}{ |c|c|c|c|c|c|c|c|c|c|c|c| }
\hline
\multirow{2}{*}{$\Omega$} & \multirow{2}{*}{$K$} & \multirow{2}{*}{$|\Omega|$} & \multirow{2}{*}{$2^R$} & \multirow{2}{*}{$\|f-\tilde{f}\|$} & \multirow{2}{*}{$\|f-\mathcal{P}_{\rT}f\|$} & \multirow{2}{*}{$\frac{\| f - \tilde{f} \|}{\| f - \cP_\rT f \|}$} & \multirow{2}{*}{$\kappa(A)$} & \multirow{2}{*}{$\frac{\sigma_{\max}(A_{4096})}{\sigma_{\min}(A)}$} & \multirow{2}{*}{$\frac{1+\delta}{\sigma_{\min}(A)}$} &  \multirow{2}{*}{$\frac{\pi}{2}\frac{1+\delta}{1-\delta}$} \\ 
 & & & & & & & & & &\\\hline
 \multirow{4}{*}{\rotatebox{90}{Jittered}} & 32 & 108 & 64 & 6.108029e-2 & 6.086270e-2 & 1.003575 &   1.550640  & 3.722720 & 4.789203 & \multirow{4}{*}{14.137167}\\
 & 64 & 215 & 128 & 3.049139e-2 & 3.046354e-2 & 1.000914 & 1.568731  & 3.840036 & 4.940129 &\\
 & 128 & 428 & 256 & 1.523943e-2 & 1.523580e-2 & 1.000238 &  1.595984  & 3.914947 & 5.036500 &\\
 & 256 & 855 & 512 & 7.618892e-3 & 7.618401e-3 &  1.000065 &   1.591625   & 4.157735 & 5.348841 &\\ \hline
 \multirow{4}{*}{\rotatebox{90}{Log}} & 32  & 350 & 64 &  6.107981e-2 & 6.086270e-2 & 1.003567 & 1.659066  &  3.415123 &  4.393487 &  \multirow{4}{*}{14.137167}\\
 & 64  & 814 &  128 & 3.049133e-2 & 3.046354e-2 & 1.000912 & 1.682514   & 3.468100&  4.461641 &\\
 & 128   & 1850 & 256 & 1.523941e-2 & 1.523580e-2 & 1.000237 &  1.694585  & 3.489929 &  4.489723 &\\
 & 256   & 4146 & 512 & 7.618890e-3 & 7.618401e-3 &  1.000064 &  1.700702  & 3.504058 &  4.507899 &\\ \hline
  \multirow{4}{*}{\rotatebox{90}{Frame}} & 32  & 76 & 64 & 6.107987e-2 & 6.086270e-2 & 1.003568 &  2.567407 & 3.445520 & \multirow{4}{*}{$\times$} &  \multirow{4}{*}{$\times$}\\
 & 64  & 144 &  128 & 3.049194e-2 & 3.046354e-2 &  1.000932 &  2.520349  & 3.318792 & &\\
 & 128   & 278 & 256 & 1.524057e-2 & 1.523580e-2 &  1.000313 &   2.621085  &  3.588619 & &\\
 & 256   & 544 & 512 & 7.618910e-3 & 7.618401e-3 & 1.000067  &   2.553133  & 3.404633 & &\\
  \hline
\end{tabular}
}
\caption{\small{The function $f(x)=\cos(6\pi x)+1/2\sin(2\pi x)$ is reconstructed by NUGS with Haar wavelets for different sampling schemes $\Omega$ and different bandwidths $K$. Jittered sampling scheme is used for $\epsilon=0.6$ and $\eta=0.1$; and log sampling scheme is used for $\delta=0.8$ and $\nu=0.4$. In the last three columns, different estimates for the reconstruction constant are computed, by using the results from \S \ref{ss:estimation} and \S \ref{ss:exact_estimate}. 
}}
\label{tab:jitt&log}
\end{table}

\begin{figure}[h] 
\footnotesize
\begin{center}
$\begin{array}{cccc}
\includegraphics[scale=0.72 ,trim=0.2cm 0.cm 0.2cm 0.2cm]{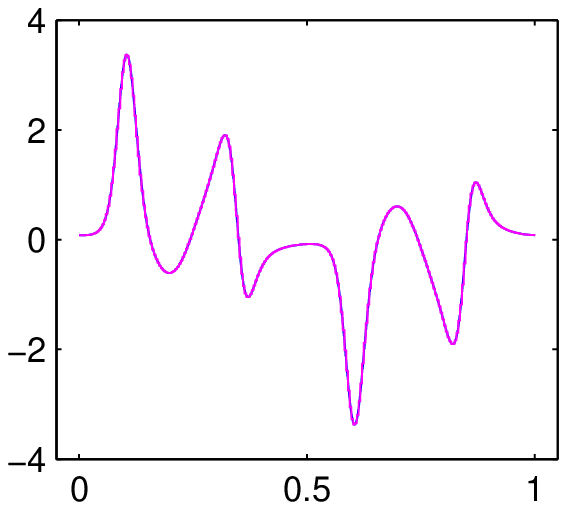} & \includegraphics[scale=0.72 ,trim=0.2cm 0.cm 0.2cm 0.2cm]{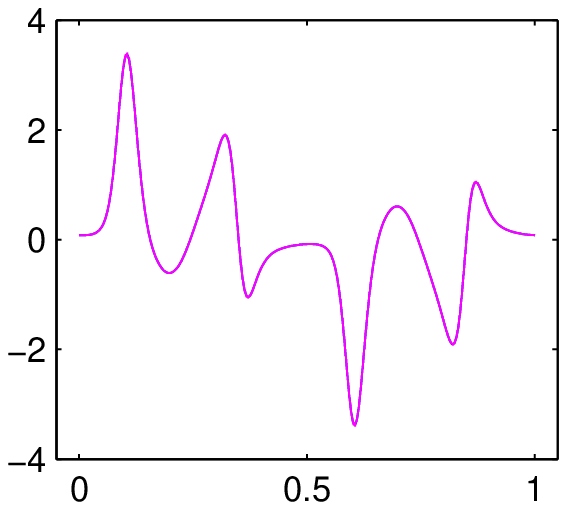} & \includegraphics[scale=0.72 ,trim=0.2cm 0.cm 0.2cm 0.2cm]{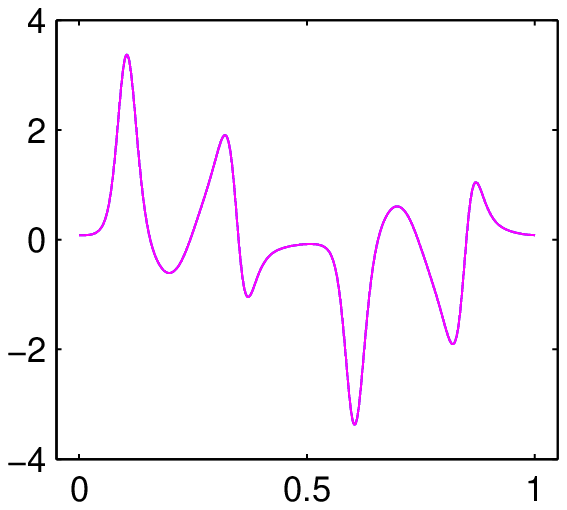}\\
\includegraphics[scale=0.72 ,trim=0.2cm 0.cm 0.2cm 0.2cm]{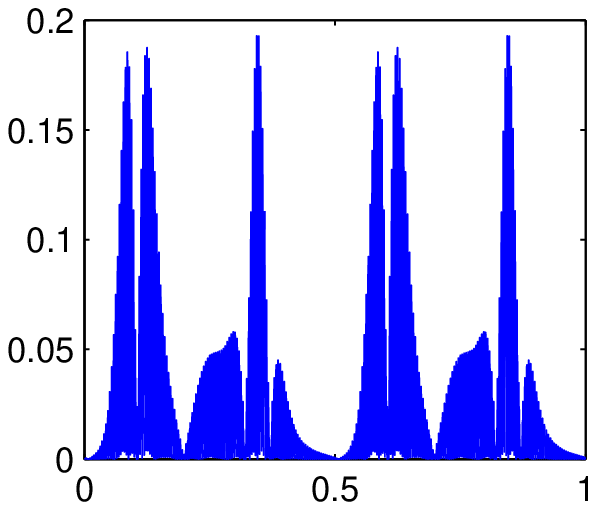} & \includegraphics[scale=0.72 ,trim=0.2cm 0.cm 0.2cm 0.2cm]{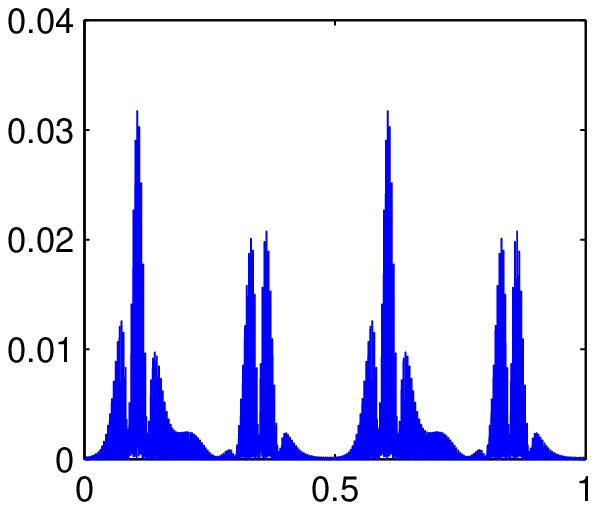} & \includegraphics[scale=0.72 ,trim=0.2cm 0.cm 0.2cm 0.2cm]{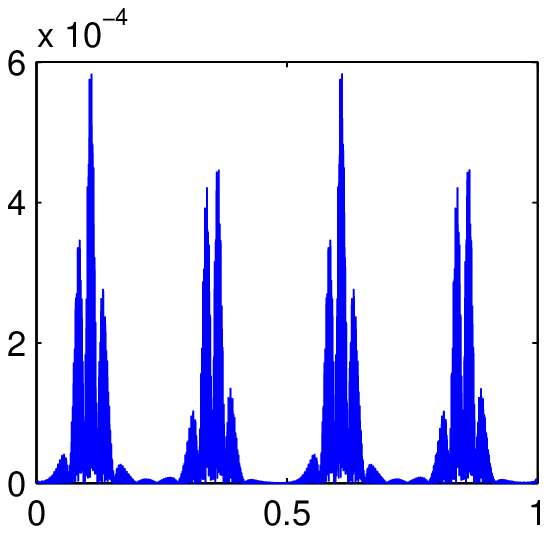}
\end{array}$\vspace{-0.3cm}
\caption{ A smooth, periodic function reconstructed by Haar, periodic DB2 and periodic DB4 wavelets, from left to right. Above is the reconstruction $\tilde{f}$ $($magenta$)$ and the original function $f$ $($blue$)$, and below is the error $|f-\tilde{f}|$. In all experiments, the same jittered sampling scheme is used, with $K=128$ and $2^R=256$.  }
\vspace{-0.2cm}
\label{fig:periodicW}
\end{center}
\end{figure}

Table \ref{tab:jitt&log} considers the case of Haar wavelet reconstructions more closely for the three different sampling schemes, and in particular, the magnitude of the reconstruction constant $C(\Omega,\rT)$.  Recall in general that $\| f - \tilde{f} \| \leq C(\Omega,\rT) \| f - \cP_{\rT} f \|$, where $\tilde{f}$ is the reconstruction.  The table suggests that this estimate is reasonably sharp.  Recall alsotechnique from \S \ref{ss:estimation} that $C(\Omega,\rT)$ can be approximated by a limiting process.  The result of this is also shown in the table.  Moreover, in the $(K,\delta)$-dense case, we see that the estimate $C(\Omega,\rT) \leq (1+\delta) / \sqrt{C_1(\Omega,\rT)}$ is also adequate (see the discussion in \S \ref{ss:estimation}).  Finally, the table also shows that the explicit bound derived in Theorem \ref{t:Haar_sharp} is also reasonably good.

\begin{figure}[h] 
\footnotesize
\begin{center}
$\begin{array}{cccc}
\includegraphics[scale=0.72 ,trim=0.2cm 0.cm 0.2cm 0.5cm]{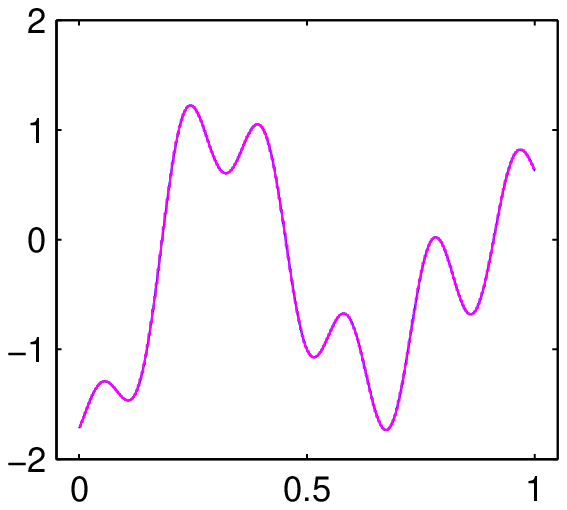} & \includegraphics[scale=0.72 ,trim=0.2cm 0.cm 0.2cm 0.5cm]{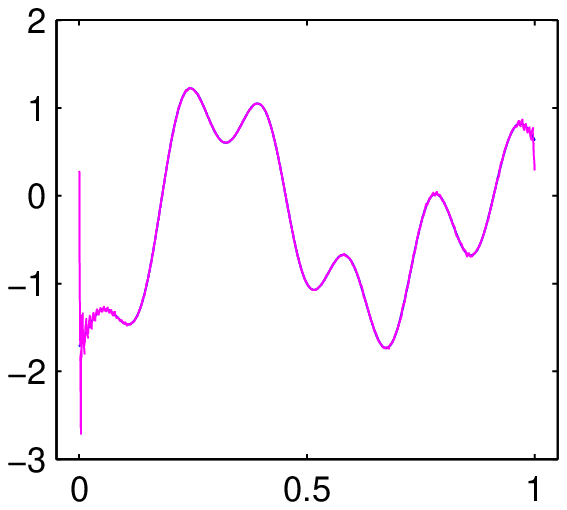} & \includegraphics[scale=0.72 ,trim=0.2cm 0.cm 0.2cm 0.5cm]{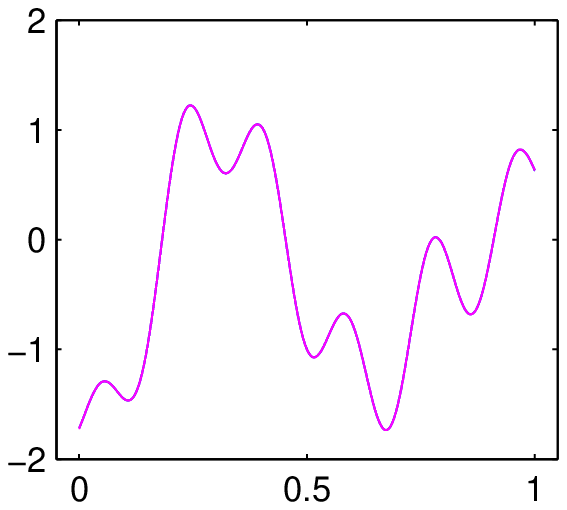}\\
\includegraphics[scale=0.72 ,trim=0.2cm 0.cm 0.2cm 0.2cm]{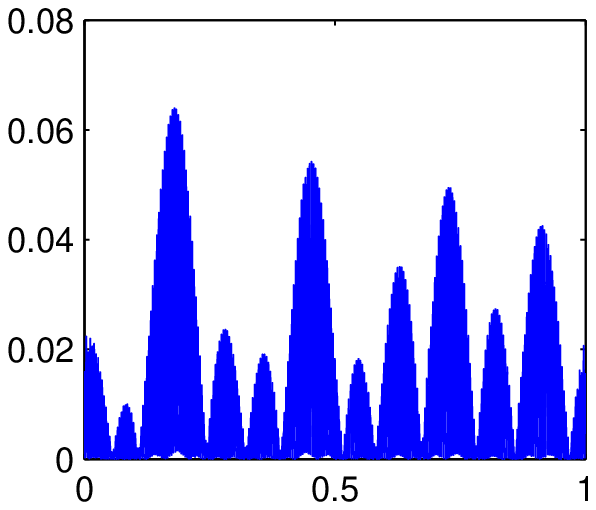} & \includegraphics[scale=0.72 ,trim=0.2cm 0.cm 0.2cm 0.2cm]{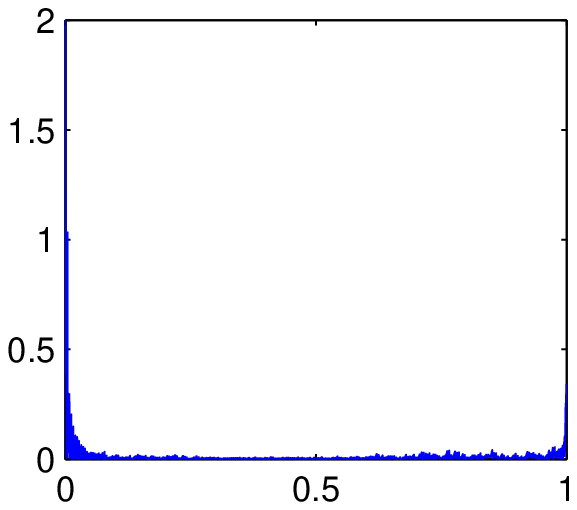} & \includegraphics[scale=0.72 ,trim=0.2cm 0.cm 0.2cm 0.2cm]{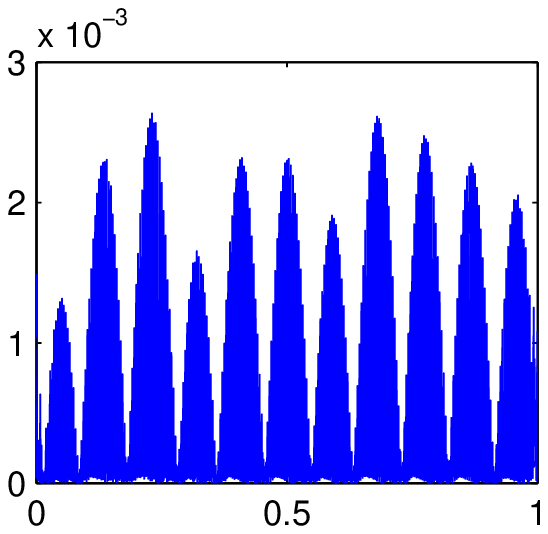}
\end{array}$\vspace{-0.3cm}
\caption{
A smooth, nonperiodic function reconstructed by Haar, periodic DB2 and boundary DB2, from left to right. Above is the reconstruction $\tilde{f}$ $($magenta$)$ and the original function $f$ $($blue$)$, and below is the error $|f-\tilde{f}|$. In all experiments, the same jittered sampling scheme is used, $K=128$ and $2^R=256$.}
\vspace{-0.2cm}
\label{fig:boundaryW}
\end{center}
\end{figure}

\begin{figure}[h]
\footnotesize
\begin{center}
$\begin{array}{cccc}
& \includegraphics[scale=0.72 ,trim=0.8cm 0.2cm 0.2cm 0.cm]{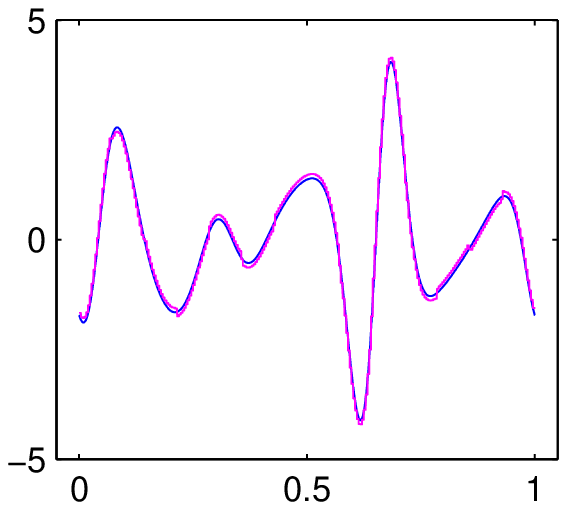} & \includegraphics[scale=0.72 ,trim=0.2cm 0.2cm 0.2cm 0.cm]{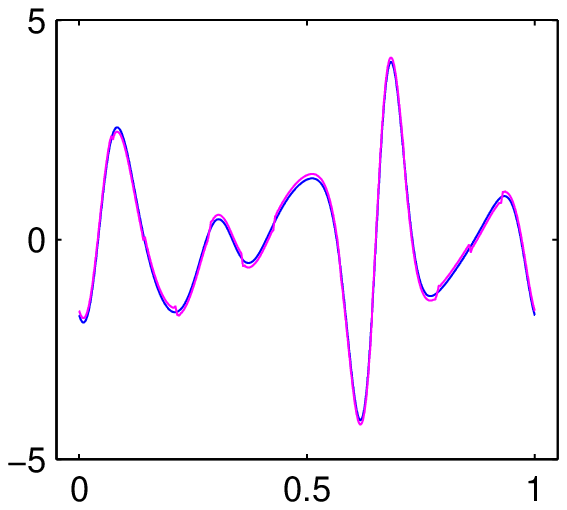} & \includegraphics[scale=0.72 ,trim=0.2cm 0.2cm 0.2cm 0.cm]{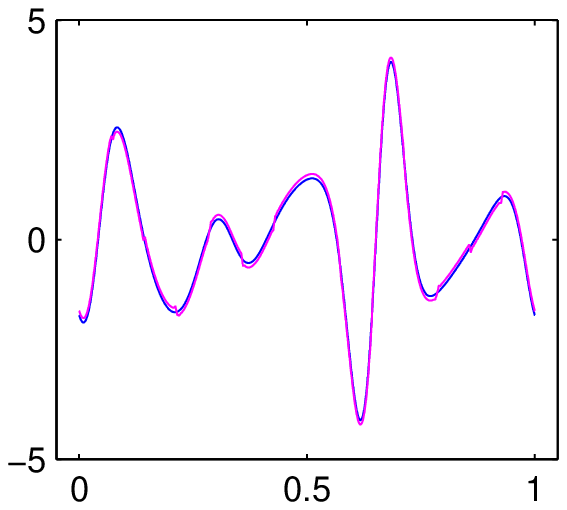}\\
\end{array}$\vspace{-0.1cm}
\caption{The function $f(x)=-\exp((\cos(6\pi x))+\sin(4\pi x))\cos(10\pi x)+\cos(4\pi x)$ $($blue$)$ and the reconstruction $F(f+\eta h)$ $($magenta$)$, where  $h(x)=\snc(14\pi(x-0.5)) \bbI_{[0,1]} /\|\snc(14\pi(x-0.5))\|$ and $\eta=0.1$. The log sampling scheme is used for $\delta=0.95$, $\nu=0.33$, $K=256$ and $N=3398$. From left to right different reconstruction basis are used for $2^R=256$: Haar, periodic DB3 and boundary DB3.  }
\vspace{-0.2cm}
\label{fig:noise}
\end{center}
\end{figure}

\begin{table}[H]
\centering
\scalebox{0.75}{
\begin{tabular}{ |c|c|c|c||c|c|c|c||c|c|c|c|c| }
\hline
$\rT$ &  $\eta$ & $\|f-F(f+\eta h)\|$ & estimate & $\rT$ &  $\eta$ & $\|f-F(f+\eta h)\|$ & estimate & $\rT$ &  $\eta$ & $\|f-F(f+\eta h)\|$ & estimate\\ \hline
\multirow{5}{*}{\rotatebox{90}{Haar}} & 0 & 4.4814e-2 & 9.4811e-2 &\multirow{5}{*}{\rotatebox{90}{DB2p}} & 0 & 3.0899e-3 & 6.5489e-3 & \multirow{5}{*}{\rotatebox{90}{DB2b}} & 0 & 4.6985e-3 & 9.6869-3\\
 &  0.05 & 6.6628e-2 & 2.0065e-1 & &  0.05  & 4.9255e-2 & 1.1259e-1 &   &  0.05  & 6.9719e-2 & 1.1521e-1\\
 &  0.1 & 1.0830e-1 & 3.0650e-1 &  &  0.1  & 9.8086e-2 & 2.1867e-1 &  &  0.1  & 1.3918e-1 & 2.2073e-1 \\
 &  0.2 & 2.0221e-1 & 5.1819e-1 & &  0.2  & 1.9609e-1 & 4.3079e-1 &  &  0.2  & 2.7826e-1 & 4.3178e-1 \\
 &  0.4 & 3.9689e-1 &  9.4158e-1 & &   0.4   & 3.9213e-1  & 8.5204e-1 &  &  0.4  & 5.5613e-1 & 8.5386e-1\\ 
  \hline
\end{tabular}
}
\caption{\small{The estimates $\tilde{C}(\Omega,\rT)\left(\|f-\cP_{\rT}f\|+\eta\|h\|\right)$ are computed for $f(x)=\cos(8\pi x)-2\sin(2\pi x)$ and $h(x)=\sin(10\pi x) \bbI_{[0,1]}/\|\sin(10\pi x)\|$, where $\tilde{C}(\Omega,\rT)=C_3(\Omega,\rT_{4096})/C_1(\Omega,\rT_{128})$ (see the Section \S \ref{ss:estimation}), and $\Omega$ is the log sampling scheme with $K=128$, $\delta=0.95$, $\nu=0.33$ and $N=1512$. The computation is done for different reconstruction spaces $\rT=\rT_{128}$ with Haar, periodic DB2 and boundary DB2 functions. } }
\label{tab:noise_new}
\end{table}

We now wish to exhibit the advantage of NUGS: namely, it allows one to reconstruct in a subspace $\rT$ that is well suited to the function to be recovered.  In Figures \ref{fig:periodicW} and \ref{fig:boundaryW} we consider the reconstruction of two functions using different wavelets.  The first function is periodic, hence we use periodic wavelets, and the second is nonperiodic, and therefore we use boundary wavelets.    Note that in all cases exactly the same set of measurements is used.  As is evident, increasing the wavelet smoothness leads to a smaller error.  This is due to the property of this approach described in Corollary \ref{c:super_NUGS}: namely, since NUGS is quasi-optimal and since it requires only a linear scaling for wavelet bases, it obtains optimal approximation rates in terms of the sampling bandwidth.

Next we consider the effect of noise.  In Table \ref{tab:noise_new} we compare the actual error in reconstructing $f$ from noisy measurements to the bound provided by $C(\Omega,\rT)$.  Note that the bound is reasonably close to the true noise value.  We also note the robustness of the reconstruction with respect to noise.  This is further illustrated in Figure \ref{fig:noise}, where we plot the reconstruction of a function $f$ from noisy measurements.  Even in the presence of large noise with $\eta  = 0.1$, we obtain a good approximation.

\subsection{Numerical comparison}
As mentioned in \S \ref{s:introduction}, two common algorithms for MRI reconstruction are gridding \cite{JacksonEtAlGridding,SedaratDCFOptimal,GelbNonuniformFourier} and iterative reconstructions \cite{FesslerFastIterativeMRI}.  We now compare these approaches with NUGS.  Recall, however, that iterative reconstruction algorithm can be interpreted as a particular instance of NUGS corresponding to a Haar wavelet basis for $\rT$ (see Remark \ref{r:it_NUGS_relation}).  We therefore continue to refer to it as such in our numerics.

\begin{figure}
\footnotesize
\begin{center}
$\begin{array}{cccc}
& \includegraphics[scale=0.72 ,trim=0.8cm 0.2cm 0.2cm 0.5cm]{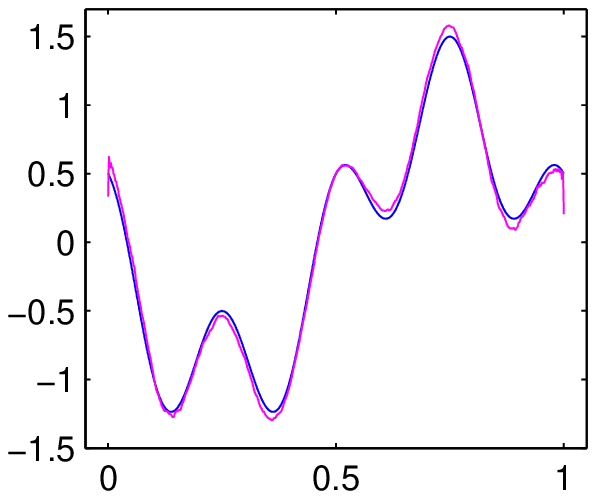} & \includegraphics[scale=0.72 ,trim=0.2cm 0.2cm 0.2cm 0.5cm]{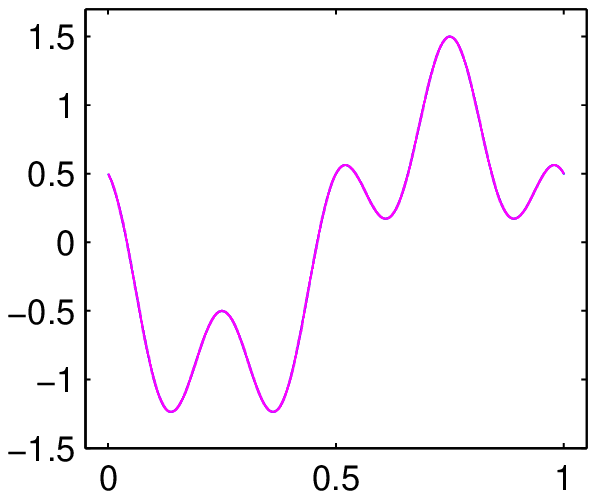} & \includegraphics[scale=0.72 ,trim=0.8cm 0.2cm 0.2cm 0.5cm]{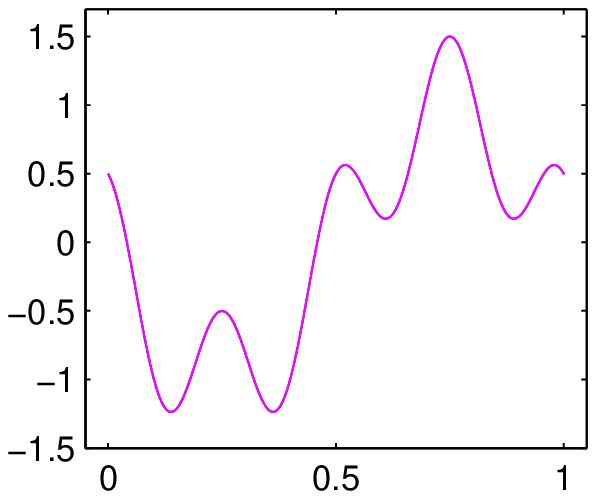}\\
& \includegraphics[scale=0.72 ,trim=0.8cm 0.2cm 0.2cm 0.5cm]{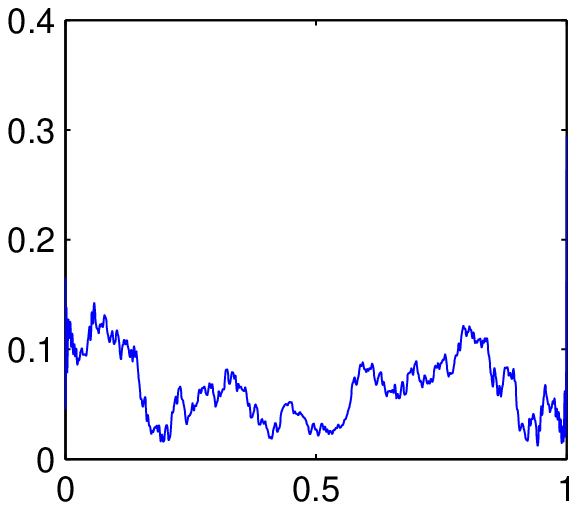} & \includegraphics[scale=0.72 ,trim=0.2cm 0.2cm 0.2cm 0.5cm]{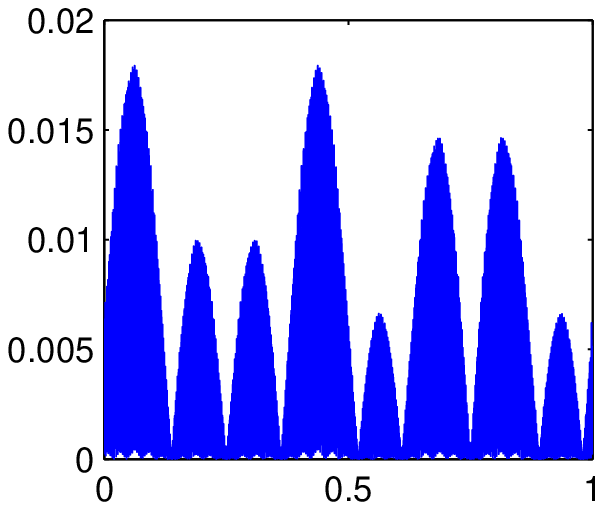} & \includegraphics[scale=0.72 ,trim=0.8cm 0.2cm 0.2cm 0.5cm]{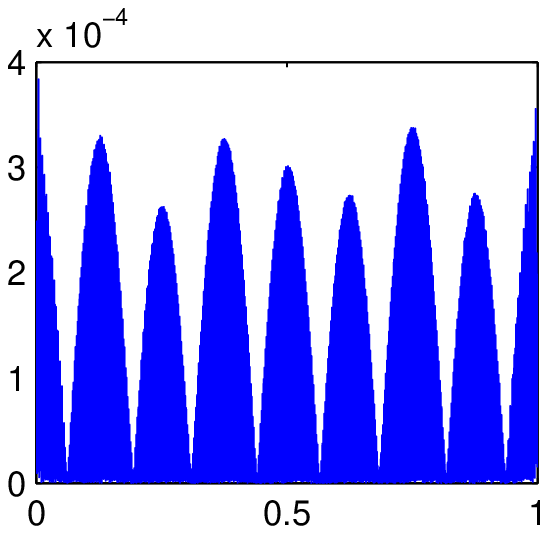}
\end{array}$\vspace{-0.1cm}
\caption{A periodic function $f(x)=1/2 \cos(8 \pi x)-\sin(2 \pi x)$ is reconstructed by gridding $($left$)$ and NUGS with Haar $($middle$)$ and DB2 $($right$)$ wavelets for $2^R=512$. The lower pictures show the error $|f-\tilde{f}|$. The jittered sampling scheme is used for  $\epsilon=0.7$, $\eta=0.14$ and $K=256$.}
\label{fig:Fourier_simple}
\end{center}
\end{figure}

Gridding is a simple technique for MRI reconstruction. It is direct, as opposed to iterative, and can be computed with a single NUFFT.  Unfortunately, this reconstruction is plagued by artefacts, even when the original function is periodic.  This is shown in the left panels of Figures \ref{fig:Fourier_simple} and \ref{fig:Four,haar,DB4}.  Alternatively, one can use the NUGS reconstruction with wavelets.  As shown in these figures, this gives a far superior reconstruction of $f$, even in the case of discontinuous functions with sharp peaks (see Figure \ref{fig:Four,haar,DB4}).  Recall also that the NUGS reconstruction can also be computed efficiently using NUFFTs (see Remark \ref{r:wavelet_fast}).  Hence, using the same measurement data, and with roughly the same computational cost, we obtain a vastly improved reconstruction.

\begin{figure}
\footnotesize
\begin{center}
$\begin{array}{cccc}
& \includegraphics[scale=0.72 ,trim=0.8cm 0.2cm 0.2cm 0.5cm]{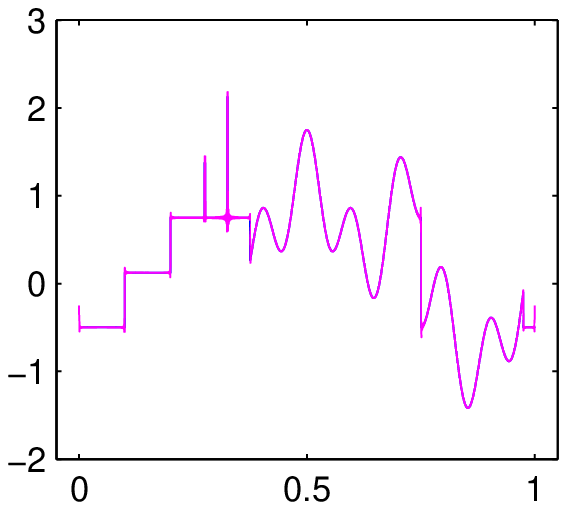} & \includegraphics[scale=0.72 ,trim=0.2cm 0.2cm 0.2cm 0.5cm]{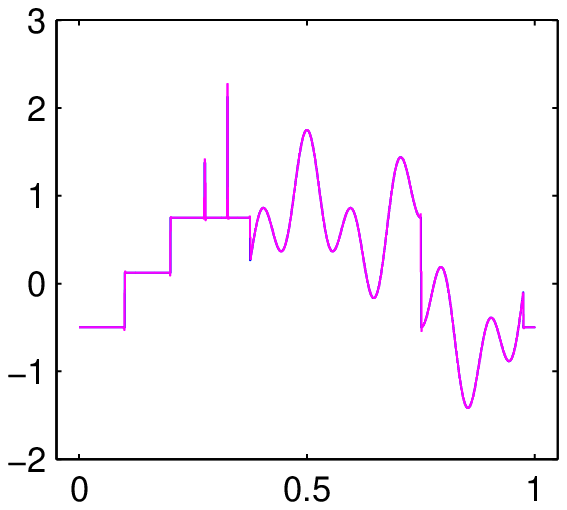} & \includegraphics[scale=0.72 ,trim=0.2cm 0.2cm 0.2cm 0.5cm]{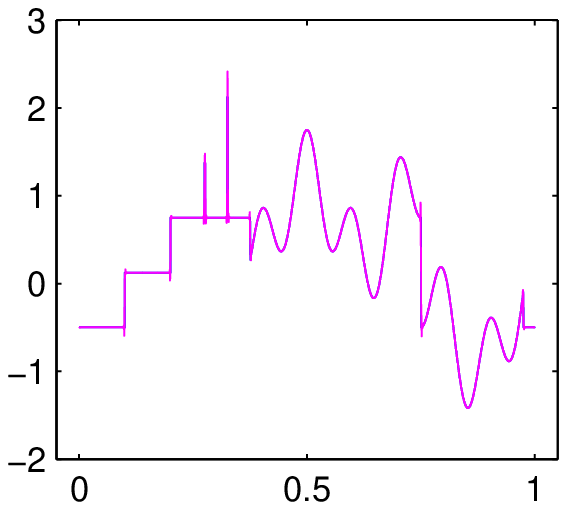}\\
& \includegraphics[scale=0.72 ,trim=0.8cm 0.2cm 0.2cm 0.5cm]{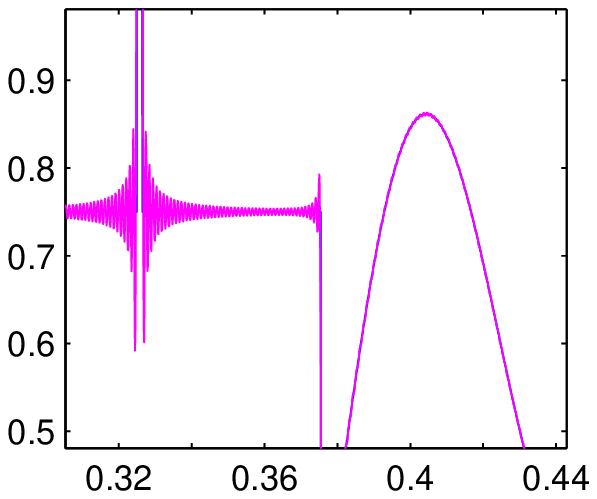} & \includegraphics[scale=0.72 ,trim=0.2cm 0.2cm 0.2cm 0.5cm]{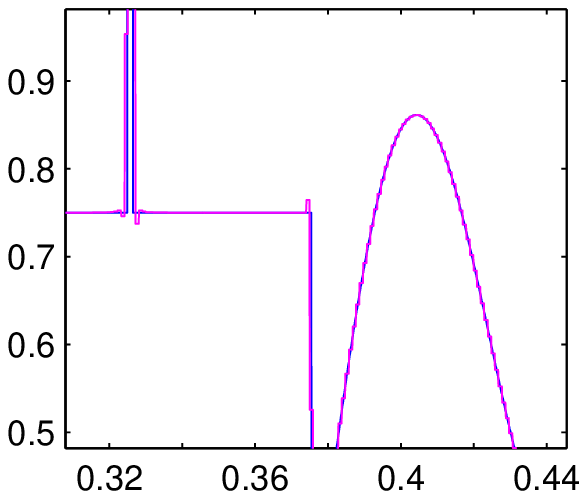} & \includegraphics[scale=0.72 ,trim=0.2cm 0.2cm 0.2cm 0.5cm]{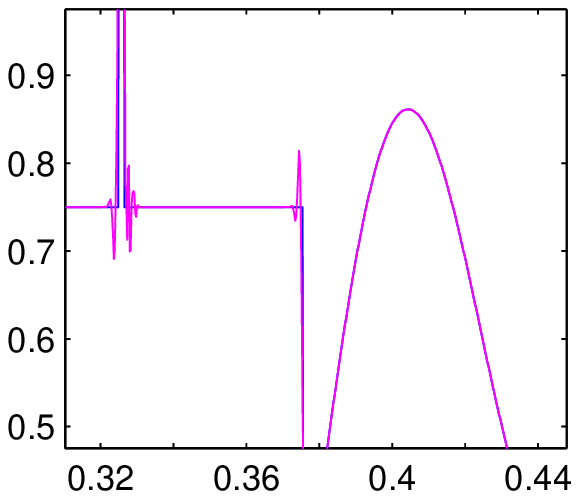}\\
\end{array}$\vspace{-0.1cm}
\caption{A discontinuous function reconstructed by gridding, and NUGS with Haar and DB4 wavelets $($from left to right$)$. The reconstruction is in magenta and original in blue. Below, a close-up is shown. The jittered sampling is used for $\epsilon=0.75$, $\eta=0.1$ and $K=2^R=1024$.
}\vspace{-0.2cm}
\label{fig:Four,haar,DB4}
\end{center}
\end{figure}

Figures \ref{fig:periodicW}--\ref{fig:Four,haar,DB4} also show the clear advantage of changing the NUGS reconstruction space $\rT$ from Haar wavelets (i.e.\ the iterative reconstructions) to higher-order wavelets.  This improvement is justified by Corollary \ref{c:super_NUGS}, following the discussion in Remark \ref{r:Why_not_DWT}.  

\rem{
The reason why NUGS obtains an improvement by changing $\rT$ is that it computes quasi-optimal approximations to the actual wavelet coefficients of $f$.  In particular, it avoids the wavelet crime \cite{StrangNguyen}.  Let $a^*$ be the vector of first $M$ coefficients in some wavelet basis. NUGS solves the least squares problem $A a \approx b$, where $A$ is the matrix of Fourier samples of wavelet basis functions and $b$ is the vector of nonuniform Fourier samples of $f$ (see \R{Abdef}).  The error estimates proved show that $\| a - a^* \| \equiv \| \tilde{f} - \cP_{\rT} f \|$ is proportional to the best approximation error $\| f - \cP_{\rT} f \|$ of $f$ in the wavelet subspace $\rT$.

As an alternative, to compute wavelet coefficients one may be tempted to construct the matrix $\tilde{A} = F W$, where $F \in \bbC^{N \times M}$ is the nonuniform discrete Fourier transform and $W \in \bbC^{M \times M}$ is the discrete wavelet transform, and solve the least squares problem $\tilde{A} a \approx b$.  Since $W$ is orthogonal, this is equivalent to solving $F c \approx b$ and then setting $a = W^T c$.  However, $c$ is a vector of pixel values of $f$, and is therefore equivalent to the solution of the iterative reconstruction algorithm (recall \S \ref{ss:exact_estimate}).  Since $W$ is orthogonal, we have $\| a - a^* \| = \| c - c^* \|$, where $c$ is the vector of exact coefficients of $f$ in the pixel basis.  Thus, the accuracy of the computed wavelet coefficients $a = W^T c$ is not determined by how well $f$ is approximated in the given wavelet basis, but how well $f$ is approximated by a piecewise constant function.  This accuracy is typically low, which means that one will not see the benefits of higher-order wavelets with this approach.  In particular, the higher approximation orders -- that is, faster decay of $\| f - \cP_{\rT} f \|$ -- offered by boundary wavelets (see Remark \ref{r:Why_not_DWT}).
}

\section{Conclusions and future work}\label{s:conclusions}
The purpose of this paper was to introduce and analyze a framework for stable reconstructions in arbitrary subspaces $\rT$ from nonuniform Fourier samples.  We have shown that this is always possible provided the samples are $(K,\delta)$-dense or arise from a Fourier frame, and provided the bandwidth $K$ or index $N$ is taken sufficiently large in relation to $\rT$.  Moreover, for the case where $\rT$ consists of wavelets, we have shown that a linear scaling of $K$ or $N$ with the dimension $M=2^R$ suffices, but that this scaling cannot be below a certain critical threshold, otherwise exponential instability necessarily occurs.

There are several topics for future work.  First, much of the one-dimensional theory developed in this paper extends to higher dimensions. In particular, following the same ideas originally due to Gr\"ochenig \cite{GrochenigIrregular}, \S 4 can be extended to higher dimensions. However, the sufficient condition of \cite{GrochenigIrregular} on the density $\delta$ for a weighted Fourier frame deteriorates linearly with dimension and thus ceases to be sharp for $d>1$.  Improving these seminal results, and in particular, establishing dimensionless estimates, is the topic of current investigations and will be presented elsewhere.  Note that in higher dimensions, it is also important to analyze other reconstruction spaces besides wavelets, such as curvelets and shearlets. 

Second, there is the question of how the reconstruction constant $C(\Omega,\rT)$ behaves for other common choices of subspace $\rT$.  In \cite{hrycakIPRM}, Hrycak \& Gr\"ochenig showed that when $\rT$ consists of polynomials of degree at most $M$, then $\ord{M^2}$ uniform Fourier samples suffice for boundedness of $C(\Omega,\rT)$ {(see also \cite{BAACHAccRecov}).  This quadratic scaling is in fact necessary, as was shown in \cite{AdcockHansenShadrinStabilityFourier}.  Similarly, when $\rT$ consists of trigonometric polynomials, it was shown in \cite{BAACHOptimality} that a linear scaling suffices whenever samples arise from a Fourier frame.  We believe both results can be extended to the $(K,\delta)$-dense case, and leave this for future work. 

A third topic for future work involves the choice of the operator $\cS$.  As discussed in Remark \ref{r:S_choice}, the theory developed in \S \ref{s:NUGS} allows for other choices of $\cS$ than that which was considered in the latter half of the paper: namely, \R{weightedFS} with weights given by \R{weights}.  It is possible that different choices, possibly depending on the subspace $\rT$, may yield improvements in the reconstruction constant.  For related work, see \cite{BergerGrochenigOblique,GelbSongFrameNFFT}.

Recall that in this paper the sampling scheme $\Omega$ is considered fixed.  This  situation arises in applications such as MRI, where $\Omega$ is often specified by physical constraints, e.g.\ magnetic gradients, noise etc.  However, in many applications, one may have substantial flexibility to design $\Omega$ so as to optimize the reconstruction quality.  That is, for a given subspace $\rT$, one seeks to design $\Omega$ as small as possible whilst keeping the reconstruction constant $C(\Omega,\rT)$ below a desired maximum value. This question is closely related to the existence of Marcinkiewicz--Zygmund inequalities (see \cite{ChuiEtAlMZ,MarzoMZ,OrtegaCerdaMZ} and references therein), which have been well-researched for certain choices of $\hat{\rT}$ (e.g.\ trigonometric polynomials, spherical harmonics,...).   On the other hand, designing good (or perhaps even optimal) sampling schemes for families of wavelet subspaces, for example, remains an open problem, but one of practical interest.

Finally, as discussed in \S \ref{s:introduction}, the eventual aim of this work is to combine the theory developed here with compressed sensing tools to allow for recovery of compressible images from relatively few nonuniform Fourier samples.  This is also work in progress.  For an extensive discussion in the case of uniform Fourier measurements we refer to \cite{AHPRBreaking}.

\appendix

\section{Proof of Theorem 3.3}

We first show that $\tilde{C} \leq C$, and in particular, that $\tilde{C} < \infty$.  By definition
\bes{
 1/\tilde{C} =\inf_{\substack{g \in \rT \\ g \neq 0}}  \frac{\| \cP_{\cS(\rT)} g \|}{\| g \|} = \inf_{\substack{g \in \rT \\ g \neq 0}} \sup_{\substack{g' \in \rT \\ \cS g' \neq 0}} \frac{|\ip{g}{\cS g'}|}{\| g \| \| \cS g' \|} .
}
Let $g \in \rT \backslash \{ 0 \}$.  If $\cS g = 0$, then $\ip{\cS g}{g}=0$ which contradicts the admissibility of $\cS$.  Hence $\cS g \neq 0$.  Therefore, we may set $g' = g$ above to get
\bes{
1/\tilde{C} \geq \inf_{\substack{g \in \rT \\ g \neq 0}} \frac{\ip{\cS g}{g}}{\| g \| \| \cS g \|}.
}
Observe that
\bes{
\| \cS g \| = \sup_{\substack{h \in \rH \\ \| h \| = 1}} \ip{\cS g}{h} \leq \sqrt{C_2} \sqrt{\ip{\cS g}{g}},
}
where the inequality follows from (3.2) and (3.4).  This now gives
\bes{
1/\tilde{C} \geq \frac{1}{\sqrt{C_2}} \inf_{\substack{g \in \rT \\ g \neq 0}} \frac{\sqrt{\ip{\cS g}{g}}}{\| g\|},
}
which, upon application of (3.3), yields $\tilde{C} \leq \sqrt{C_2/C_1} = C$ as required.

To prove the remainder of the theorem, we shall used the techniques of [5] based on the geometric notions of subspace angles and oblique projections.  Let $\rU = \rT$ and $\rV = (\cS(\rT))^{\perp}$. Note that $1/\tilde{C} = \cos (\theta_{\rU \rV^{\perp}})$ is cosine of the subspace angle between $\rU$ and $\rV^{\perp}$ defined by
\bes{
\cos (\theta_{\rU \rV^{\perp}}) = \inf_{\substack{u \in \rU \\ \| u \|=1}} \| \cP_{V^{\perp}} u \|.
}
Since $\tilde{C}<\infty$, the subspaces $\rU$ and $\rV$ satisfy the so-called subspace condition $\cos (\theta_{\rU \rV^{\perp}}) > 0$.  Thus [5,\ Cor.\ 3.5] gives
\bes{
\| \cW_{\rU \rV} f \| \leq \tilde{C} \| f \|,\quad \forall f \in \rH_0,
}
and
\bes{
\| f - \cW_{\rU \rV} f \| \leq \tilde{C} \| f - \cP_{\rU} f \|,\quad \forall f \in \rH_0,
}
where $\rH_0 = \rU \oplus \rV$ and $\cW_{\rU \rV} : \rH_0 \rightarrow \rU$ is the projection with range $\rU$ and kernel $\rV$.  

Hence, to establish (3.9) it remains to show the following: (i) $\rH_0 = \rH$ and (ii) $\tilde{f} = \cW_{\rU \rV} f$, $\forall f \in \rH$.  For (i), we note that $\rH_0 = \rH$ provided $\dim(\cS(\rT)) = \dim(\rT)$  [5,\ Lem.\ 3.10].  However, if not then there exists a nonzero $g \in \rT$ such that $\cS(g) = 0$.  As previously observed, this implies that $g = 0$; a contradiction.

For (ii), we first note that
\bes{
\ip{\cW_{\rU \rV} f }{\cS g} = \ip{f}{\cS g},\quad \forall g\in \rT.
}
Since $\cS$ is self-adjoint, it follows that $\cW_{\rU \rV} f$ satisfies the same conditions (3.7) as $\tilde{f}$.  Thus, it remains only to show that $\tilde{f}$ is unique.  However, if not then we find that there is a nonzero $g \in \rT \cap \cS(\rT)^{\perp} = \rU \cap \rV$.  But then $\cos (\theta_{\rU \rV^\perp}) = 0$, and this contradicts the fact that $\rU$ and $\rV$ obey the subspace condition.

\section{Construction of periodic, folded and boundary wavelets}

Following [48], 
we consider three standard constructions -- periodic, folded and boundary wavelets.

\subsubsection*{Periodic wavelets}\label{sss:periodic}
Suppose that $\{ \psi_{j,k} \}_{j,k \in \bbZ}$ is a wavelet basis of $\rL^2(\bbR)$ associated to an MRA with scaling function $\phi$.  Define the periodizing operation
\bes{
f(x) \mapsto f^{\mathrm{per}}(x) = \sum_{k \in \bbZ} f(x+k),
}
and let $\psi^{\mathrm{per}}_{j,k}$ and $\phi^{\mathrm{per}}_{j,k}$ be the corresponding periodic wavelets and scaling functions.  Define the periodized MRA spaces
\bes{
V^{\per}_j = \spn \left \{ \phi^{\per}_{j,k} :\ k=0,\ldots,2^{j}-1 \right \} ,\quad W^{\per}_j = \spn \left \{ \psi^{\per}_{j,k} : \ k=0,\ldots,2^{j}-1 \right \}.
}
Note that the maximal index $k$ is finite, since $\phi^{\per}_{j,k+2^j} = \phi^{\per}_{j,k}$ and likewise for $\psi^{\per}_{j,k}$.  

Now let $J \in \bbN_0$ be given.  Then 
\bes{
\rL^2(0,1) = \overline{V^{\per}_{J}  \oplus W^\per_J \oplus W^{\per}_{J+1} \oplus \cdots}.
}
We now introduce the finite-dimensional reconstruction space $\rT$:
\bes{
\rT = V^{\per}_J \oplus W^\per_J \oplus W^{\per}_{J+1} \oplus \cdots \oplus W^{\per}_{R-1}.
}
This is the space (5.2) for $\text{type}=\per$.
Since the original wavelets have an MRA, we also have that
\bes{
\rT = V^{\per}_R = \spn \left \{ \phi^{\per}_{R,k} : k=0,\ldots,2^R-1 \right \}.
}
Without loss of generality, we now suppose that $\supp(\phi) \subseteq [-p+1,p]$ for $p\in \bbN$.  Note the following: if $\supp(f) \subseteq [0,1]$ then $f(x) = f^{\per}(x)$ for $x \in [0,1]$.  In particular, since
\bes{
\supp (\phi_{R,k}) = [(k-p+1)/ 2^R , (k+p)/2^R ],
}
we have that $\phi^{\per}_{R,k}(x) = \phi_{R,k}(x)$, $x \in [0,1]$, whenever $k=p,\ldots,2^R-p-1$.  Hence, it is possible to decompose the space $\rT$ into
\bes{
\rT = \rT^{\lt} \oplus \rT^i \oplus \rT^{\rt},
}
where 
\bes{
\rT^{i} = \spn \left \{ \phi_{R,k} : k=p,\ldots,2^R-p-1 \right \},
}
contains interior  scaling functions with support in $(0,1)$ and 
\eas{
\rT^{\lt} &= \spn \left \{ \phi^{\per}_{R,k} \bbI_{[0,1]} : k=0,\ldots,p-1 \right \}, \\
\rT^{\rt} &= \spn \left \{  \phi^{\per}_{R,k}  \bbI_{[0,1]}  : k = 2^R-p,\ldots,2^R-1 \right \},
}
contain the periodized scaling functions.

\subsubsection*{Folded wavelets}\label{sss:folded}

Folded wavelets are defined via the folding operation
\bes{
f(x) \mapsto f^{\mathrm{fold}}(x) = \sum_{k \in \bbZ} f(x-2k) + \sum_{k \in \bbZ} f(2k-x).
}
In this case, one obtains biorthogonal bases of wavelets for $\rH$.  Note that in this case we have
\bes{
V^{\fold}_j = \spn \left \{ \phi^{\fold}_{j,k} :\ k=0,\ldots,2^{j}-\iota \right \} ,\quad W^{\fold}_j = \spn \left \{ \psi^{\fold}_{j,k} : \ k=0,\ldots,2^{j}-1 \right \},
}
where $\iota$ takes value $0$ if the wavelets are symmetric about $x=1/2$ and $1$ if they are antisymmetric.  Much as before, we define the finite-dimensional reconstruction space
\bes{
\rT = V^{\fold}_J \oplus W^\fold_J \oplus W^{\fold}_{J+1} \oplus \cdots \oplus W^{\fold}_{R-1},
}
and note that
\bes{
\rT = V^{\fold}_{R} =\spn \left \{ \phi^{\fold}_{R,k} : k=0,\ldots,2^R-\iota\right \}.
}
As for periodic wavelets, $\rT$ can be decomposed into three subspaces containing interior and boundary wavelets, i.e. $\rT=\rT^{\lt} \oplus \rT^i \oplus \rT^{\rt}$.

\subsubsection*{Boundary wavelets}\label{sss:boundary}
In this case, a new basis on $[0,1]$ is defined as follows. Set
\bes{
\phi^{\ivl}_{j,k}(x) = \left \{ \begin{array}{ll} 2^{j/2} \phi (2^j x - k) & p \leq k < 2^j - p 
\\ 
2^{j/2} \phi^{\lt}_{k} (2^j x) & 0 \leq k <p 
\\
2^{j/2} \phi^{\rt}_{2^j-k-1}(2^j(x-1)) & 2^j - p \leq k < 2^j,
\end{array} \right .
}
and similarly for the wavelet functions $\psi^{\ivl}_{j,k}$.  Here the functions $\phi^{\lt}_k$ and $\phi^{\rt}_k$ are particular boundary scaling functions.  
See [21] 
for details. We may now define an MRA
\bes{
V^{\ivl}_{j} = \spn \left \{ \phi^{\ivl}_{j,k} :\ k=0,\ldots,2^j-1 \right \},\quad W^{\ivl}_{j} = \spn \left \{  \phi^{\ivl}_{j,k} :\ k=0,\ldots,2^j-1 \right \},
}
which, for $J \geq \log_2(2p)$ gives the reconstruction space
\bes{
\rT = V^{\ivl}_{J} \oplus W^{\ivl}_J \oplus \cdots \oplus W^{\ivl}_{R-1} = V^{\ivl}_{R}.
}
Note that, as before, we may decompose  $\rT$ as $\rT^{\lt} \oplus \rT^i \oplus \rT^{\rt}$.

\section{Proof of Proposition 5.12}

First, we require the following lemma:

\lem{
\label{l:augment}
Let $\{ \omega_n \}_{n \in \bbZ}$ be an increasing sequence of separated points with minimal separation $\eta = \inf_{n \in \bbZ} \{ \omega_{n+1} - \omega_n \} > 0 $.  Then there exists a set of points $\{ \tilde{\omega}_n \}_{n \in \bbZ}$ with minimal separation at least $\eta/2$ such that $\{ \omega_n \}_{n \in \bbZ} \subseteq \{ \tilde{\omega}_n \}_{n \in \bbZ}$ and
\bes{
\sup_{n \in \bbZ} \{  \tilde{\omega}_{n+1} - \tilde{\omega}_{n} \} \leq \eta.
}
}
\prf{
Let $n \in \bbZ$.  If $\omega_{n+1} - \omega_n = \eta$ then we do nothing.  Otherwise, let $k \in \bbN$ be the smallest integer such that $\omega_{n+1} - \omega_{n} \leq (k+1) \eta$.  Introduce the new points
\bes{
\omega_{n} + r \eta,\quad r=1,\ldots,k-1,
}
as well as
\bes{
\frac12 \left ( \omega_{n} + (k-1) \eta + \omega_{n+1} \right ).
}
These new points are at least $\eta/2$ separated, and have maximal separation at most $\eta$.
}

\prf{[Proof of Proposition 5.12.]
Recall from Theorem 4.7 that any sequence $\{ \omega_n \}_{n \in \bbZ}$ that gives a frame is necessarily relatively separated, i.e.\ it is a finite union of separated sequences.  Since we wish to obtain an upper bound for
\bes{
\sum_{|n| > N} | \hat{f}(\omega_n) |^2,
}
for any $f \in \rT$, we may therefore assume without loss of generality that $\{ \omega_n \}_{n \in \bbZ}$ is a separated sequence with separation $\eta$.  Moreover, after an application of Lemma \ref{l:augment}, we may assume without loss of generality that $\{ \omega_n \}_{n \in \bbZ}$ is $\eta / 2$ separated with maximal spacing at most $\eta$.

As in the proof of Proposition 5.10, let $f = \sum^{M_2}_{k=M_1} a_k \sqrt{M} \phi(M \cdot - k ) \in \rT$ and write $\tilde{\Psi}(x) =\linebreak \sum^{M_2}_{k=M_1} a_k \E^{-2 \pi \I k x}$ so that
\be{
\label{detroit}
\hat{f}(\omega) = \frac{1}{\sqrt{M}} \hat{\phi} \left ( \frac{\omega}{M} \right ) \tilde \Psi \left ( \frac{\omega}{M} \right ).
}
Let
\be{\label{sym_trig}
\Psi(x) = \E^{2\pi \I M_3 x} \tilde{\Psi}(x) = \sum^{M_2-M_3}_{k=M_1-M_3} a_{k+M_3} \E^{-2 \pi \I k x},\qquad M_3  =  \left \lceil \frac{M_1+M_2}{2} \right \rceil,
}
so that $| \Psi(x) | = | \tilde{\Psi}(x) |$.  By \R{detroit} we also have $|\hat{f}(\omega) | =\frac{1}{\sqrt{M}} | \hat{\phi}(\omega/M) | | \tilde{\Psi}(\omega/M ) |  =\frac{1}{\sqrt{M}} | \hat{\phi}(\omega/M) | | \Psi(\omega/M ) |$, and therefore
\eas{
\sum_{n > N} | \hat{f}(\omega_n) |^2 &\leq \frac{1}{M} \sum_{n > N}  | \hat{\phi}(\omega_n/M) |^2 | \Psi(\omega_n/M ) |^2
\\
& \leq \frac{1}{M} \sum^{\infty}_{l=0} \sup_{\omega \in I_l} | \hat{\phi}(\omega/M) |^2 \sum_{\substack{n : \omega_n \in I_l}} | \Psi(\omega_n / M) |^2,
}
where $I_l = [\omega_{N} + l M , \omega_N + (l+1) M )$.  Since $\{ \omega_n \}_{n \in \bbZ}$ is separated and increasing, we must have that $\omega_N \gtrsim N$ as $N \rightarrow \infty$.  In particular $\omega_N > 0$ for sufficiently large $N$.  By the assumption on $\phi$, we therefore obtain
\bes{
\sum_{n > N} | \hat{f}(\omega_n) |^2 \lesssim  M^{2\alpha - 1} \sum^{\infty}_{l=0} \left ( \omega_N + 2 l M \right )^{-2 \alpha}  \sum_{\substack{n : \omega_n \in I_l}} | \Psi(\omega_n / M) |^2.
}
We now claim that the result follows, provided
\be{
\label{claim}
\sum_{\substack{n : \omega_n \in I_l}} | \Psi(\omega_n / M) |^2 \leq c M \| \Psi \|^2,\quad \forall l =0,1,2,\ldots.
}
We shall prove that \R{claim} holds in a moment.  First, however, let us show how \R{claim} implies the result.  Substituting this bound into the previous expression gives
\eas{
\sum_{n > N} | \hat{f}(\omega_n) |^2 \lesssim M^{2 \alpha} \sum^{\infty}_{l=0} (\omega_N + 2 l M )^{-2 \alpha} \| \Psi \|^2 \lesssim \left ( \frac{\omega_N}{M} \right )^{1-2 \alpha} \| \Psi \|^2.
}
Similarly, we also get
\bes{
\sum_{n < -N} | \hat{f}(\omega_n) |^2 \lesssim \left ( \frac{|\omega_{-N}|}{M} \right )^{1-2 \alpha} \| \Psi \|^2.
}
An application of (5.9) now gives
\bes{
\tilde{E}(\rT,N)^2 \lesssim \frac{1}{d_1} \left ( \frac{\min\{ \omega_N, | \omega_{-N} | \}}{M} \right )^{1-2 \alpha}.
}
Since $\omega_N, | \omega_{-N} | \gtrsim N$ as $N \rightarrow \infty$, the result now follows.

It remains to establish \R{claim}.  Write $\{ \omega_n / M : \omega_n \in I_l \} = \{ x_1,\ldots,x_L \}$ where
\bes{
\omega_{N}/M + l \leq x_1 < x_2 < \ldots < x_L \leq \omega_{N}/M + l+1,
}
and set $x_0 = x_1$ and $x_{L+1} = x_L$.  Note that $\eta / (2M) \leq x_{n+1} - x_n \leq \eta / M$.  Therefore
\bes{
\sum_{\substack{n : \omega_n \in I_l}} | \Psi(\omega_n / M) |^2  = \sum^{L}_{n=1} | \Psi(x_n) |^2 \leq \frac{2 M}{\eta} \sum^{L}_{n=1} \mu_n | \Psi(x_n) |^2,
}  
where $\mu_n = \frac12(x_{n+1}-x_{n-1})$. Hence, by Lemma 5.11 we have
\bes{
\sum_{\substack{n : \omega_n \in I_l}} | \Psi(\omega_n / M) |^2  \leq \frac{2 M}{\eta}\left [ \| \Psi \|_{[a,b]} + \frac{\eta}{M \pi} \| \Psi' \|_{[a,b]} \right ]^2,
}
where $a = \frac12(x_1+x_0) = x_1$ and $b = \frac12(x_{L+1}+x_L) = x_L$.  Note that $|b-a| \leq 1$.  Hence since $\Psi$ is periodic, we get 
\eas{
\sum_{\substack{n : \omega_n \in I_l}} | \Psi(\omega_n / M) |^2  \leq \frac{2 M}{\eta} \left [ \| \Psi \| + \frac{\eta}{M \pi} \| \Psi' \| \right ]^2.
}
To prove the result, we only need to show that $\| \Psi' \| \leq M \pi \| \Psi \|$.  Since $\Psi$ is a trigonometric polynomial given by \R{sym_trig}, we have
\bes{
\| \Psi' \| \leq 2 \max \left \{ M_2 - M_3 , M_3 - M_1 \right \} \pi \| \Psi \|.
}  
Thus it remains to show that $M_2 - M_3 , M_3 - M_1 \leq M/2$.  Since $\rT \subseteq \rH$ by assumption, the function $\phi$ must have compact support.  Let $\supp(\phi) \subseteq [a,b]$.  Then we must also have that $-a \leq M_1 \leq M_2 \leq M-b$.  In particular, $M_2 - M_1 \leq M - (b-a) < M$.  Therefore
\bes{
M_2 - M_3 \leq M_2 - \frac{M_1+M_2}{2} <\frac{M}{2},\quad M_3 - M_1 \leq \frac{M_1+M_2}{2} + 1 - M_1 \leq \frac{M}{2} + 1 - \frac{b-a}{2}.
}
Since $M_3 - M_1 \in \bbN$ and $b - a >0$ we obtain the result.
}

\section*{Acknowledgements}
The authors would like to thank Akram Aldroubi, Anne Gelb, Karlheinz Gr\"ochenig, Rodrigo Platte and Yang Wang for useful discussions.  BA acknowledges support from the NSF DMS grant 1318894.  MG acknowledges support from the UK Engineering and Physical Sciences Research Council (EPSRC) grant EP/H023348/1 for the University of Cambridge Centre for Doctoral Training, the Cambridge Centre for Analysis.  AH acknowledges support from a Royal Society University Research Fellowship as well as the EPSRC grant EP/L003457/1.

\bibliographystyle{abbrv}
\small
\bibliography{Nonuniform1DRefs}

\begin{thebibliography}{10}

\bibitem{BAACHGSCS}
B.~Adcock and A.~C. Hansen.
\newblock {Generalized Sampling and Infinite-Dimensional Compressed Sensing}.
\newblock {\em Technical report NA2011/02, DAMTP, University of Cambridge},
  2011.

\bibitem{BAACHShannon}
B.~Adcock and A.~C. Hansen.
\newblock A generalized sampling theorem for stable reconstructions in
  arbitrary bases.
\newblock {\em J. Fourier Anal. Appl.}, 18(4):685--716, 2012.

\bibitem{BAACHAccRecov}
B.~Adcock and A.~C. Hansen.
\newblock {Stable reconstructions in {H}ilbert spaces and the resolution of the
  {G}ibbs phenomenon}.
\newblock {\em Appl. Comput. Harmon. Anal.}, 32(3):357--388, 2012.

\bibitem{AHKM2DWavelets}
B.~Adcock, A.~C. Hansen, G.~Kutyniok, and J.~Ma.
\newblock Linear stable sampling rate: Optimality of {2D} wavelet
  reconstructions from {F}ourier measurements.
\newblock {\em arXiv:1403.0172}, 2014.

\bibitem{BAACHOptimality}
B.~Adcock, A.~C. Hansen, and C.~Poon.
\newblock Beyond consistent reconstructions: optimality and sharp bounds for
  generalized sampling, and application to the uniform resampling problem.
\newblock {\em SIAM J. Math. Anal.}, 45(5):3114--3131, 2013.

\bibitem{AHPWavelet}
B.~Adcock, A.~C. Hansen, and C.~Poon.
\newblock On optimal wavelet reconstructions from {F}ourier samples: linearity
  and universality of the stable sampling rate.
\newblock {\em Appl. Comput. Harmon. Anal. (to appear)}, 2013.

\bibitem{AHPRBreaking}
B.~Adcock, A.~C. Hansen, C.~Poon, and B.~Roman.
\newblock Breaking the coherence barrier: A new theory for compressed sensing.
\newblock {\em arXiv:1302.0561}, 2014.

\bibitem{AHRT2013AEIP}
B.~Adcock, A.~C. Hansen, B.~Roman, and G.~Teschke.
\newblock Generalized sampling: stable reconstructions, inverse problems and
  compressed sensing over the continuum.
\newblock {\em Advances in Imaging and Electron Physics (to appear)}, 2013.

\bibitem{AdcockHansenShadrinStabilityFourier}
B.~Adcock, A.~C. Hansen, and A.~Shadrin.
\newblock A stability barrier for reconstructions from {F}ourier samples.
\newblock {\em SIAM J. Numer. Anal.}, 52(1):125--139, 2014.

\bibitem{AhnEtAlSpiralNMR}
C.~B. Ahn, J.~H. Kim, and Z.~H. Cho.
\newblock {High-speed spiral-scan echo planar {NMR} imaging-{I}}.
\newblock {\em IEEE Trans. Med. Imaging}, 5(1):2--7, 1986.

\bibitem{AldroubiAverageSamp}
A.~Aldroubi.
\newblock Non-uniform weighted average sampling and reconstruction in
  shift-invariant and wavelet spaces.
\newblock {\em Appl. Comput. Harmon. Anal.}, 13:151--161, 2002.

\bibitem{AldroubiGrochenigSIREV}
A.~Aldroubi and K.~Gr{\"o}chenig.
\newblock {Nonuniform sampling and reconstruction in shift-invariant spaces}.
\newblock {\em SIAM Rev.}, 43:585--620, 2001.

\bibitem{AldroubiConvolutionAvg}
A.~Aldroubi, Q.~Sun, and W.-S. Tang.
\newblock Convolution, average sampling and a {C}alderon resolution of the
  identity for shift-invariant spaces.
\newblock {\em J. Fourier Anal. Appl.}, 22:215--244, 2005.

\bibitem{MRI1D2}
O.~Beckonert, M.~Coen, H.~Keun, Y.~Wang, T.~Ebbels, E.~Holmes, J.~Lindon, and
  J.~Nicholson.
\newblock High-resolution magic-angle-spinning nmr spectroscopy for metabolic
  profiling of intact tissues.
\newblock {\em Nature Protocols}, 5(6):1019 -- 1032, 2010.

\bibitem{BenedettoIrregular}
J.~J. Benedetto.
\newblock {Irregular sampling and frames}.
\newblock In J.~J. Benedetto and M.~Frazier, editors, {\em {Wavelets:
  Mathematics and Applications}}. Boca Raton, FL: CRC, 1994.

\bibitem{BenedettoSpiral}
J.~J. Benedetto and H.~C. Wu.
\newblock {Non-uniform sampling and spiral {MRI} reconstruction}.
\newblock {\em Proc. SPIE}, 4119:130--141, 2000.

\bibitem{BergerGrochenigOblique}
P.~Berger and K.~Gr\"{o}chenig.
\newblock Sampling and reconstruction in different subspaces by using oblique
  projections.
\newblock {\em arXiv:1312.1717}, 2013.

\bibitem{ChristensenFramesAMS}
O.~Christensen.
\newblock Frames, {R}iesz bases, and discrete {G}abor/wavelet expansions.
\newblock {\em Bull. Amer. Math. Soc}, 38(3):273--291, 2001.

\bibitem{ChuiEtAlMZ}
C.~K. Chui and L.~Zhong.
\newblock {Polynomial interpolation and Marcinkiewicz--Zygmund inequalities on
  the unit circle}.
\newblock {\em J. Math. Anal. Appl.}, 233(1):387--405, 1999.

\bibitem{CDF}
A.~Cohen, I.~Daubechies, and J.~Feauveau.
\newblock Biorthogonal bases of compactly supported wavelets.
\newblock {\em Comm. Pure Appl. Math.}, 45(5):485--560, 1992.

\bibitem{CDV}
A.~Cohen, I.~Daubechies, and P.~Vial.
\newblock {Wavelets on the Interval and Fast Wavelet Transforms}.
\newblock {\em Appl. Comput. Harmon. Anal.}, 1(1):54 -- 81, 1993.

\bibitem{DelattreEtAlSprial}
B.~M.~A. Delattre, R.~M. Heidemann, L.~A. Crowe, Vall{\'e}e, and J.-N.
  Hyacinthe.
\newblock {Spiral demystified}.
\newblock {\em Magn. Reson. Imaging}, 28(862--881), 2010.

\bibitem{MRI1D1}
T.~Diki\'{c}, S.~J.~F. Erich, W.~Ming, H.~P. Huinink, P.~C. Thüne, R.~A. T.~M.
  van Benthem, and G.~de~With.
\newblock Fluorine depth profiling by high-resolution 1d magnetic resonance
  imaging.
\newblock {\em Polymer}, 48(14):4063 -- 4067, 2007.

\bibitem{NMR3}
K.~Eberhardt, C.~Degen, A.~Hunkeler, and B.~Meier.
\newblock One- and two-dimensional nmr spectroscopy with a magnetic-resonance
  force microscope.
\newblock {\em Angewandte Chemie International Edition}, 47(46):8961--8963,
  2008.

\bibitem{eldar2003sampling}
Y.~C. Eldar.
\newblock Sampling without input constraints: Consistent reconstruction in
  arbitrary spaces.
\newblock In A.~I. Zayed and J.~J. Benedetto, editors, {\em Sampling, Wavelets
  and Tomography}, pages 33--60. Boston, MA: Birkh{\"a}user, 2004.

\bibitem{eldar2005general}
Y.~C. Eldar and T.~Werther.
\newblock General framework for consistent sampling in {H}ilbert spaces.
\newblock {\em Int. J. Wavelets Multiresolut. Inf. Process.}, 3(3):347, 2005.

\bibitem{FeichtingerGrochenigIrregular}
H.~G. Feichtinger and K.~Gr{\"o}chenig.
\newblock {Theory and practice of irregular sampling}.
\newblock In J.~J. Benedetto and M.~Frazier, editors, {\em {Wavelets:
  Mathematics and Applications}}, pages 305--363. Boca Raton, FL: CRC, 1994.

\bibitem{FeichtingerEtAlEfficientNonuniform}
H.~G. Feichtinger, K.~Gr{\"o}chenig, and T.~Strohmer.
\newblock {Efficient numerical methods in nonuniform sampling theory}.
\newblock {\em Numer. Math.}, 69:423--440, 1995.

\bibitem{GelbSongFrameNFFT}
A.~Gelb and G.~Song.
\newblock {A frame theoretic approach to the Non-Uniform Fast Fourier
  Transform}.
\newblock {\em SIAM J. Numer. Anal. (to appear)}, 2014.

\bibitem{GrochenigIrregular}
K.~Gr{\"o}chenig.
\newblock {Reconstruction algorithms in irregular sampling}.
\newblock {\em Math. Comp.}, 59:181--194, 1992.

\bibitem{GrochenigLAA}
K.~Gr\"{o}chenig.
\newblock A discrete theory of irregular sampling.
\newblock {\em Linear Algebra Appl.}, 193:129--150, 1993.

\bibitem{GrochenigIrregularExpType}
K.~Gr\"{o}chenig.
\newblock Irregular sampling, {T}oeplitz matrices, and the approximation of
  entire functions of exponential type.
\newblock {\em Math. Comp.}, 68(226):749--765, 1999.

\bibitem{GrochenigModernSamplingBook}
K.~Gr\"{o}chenig.
\newblock Non-uniform sampling in higher dimensions: from trigonometric
  polynomials to bandlimited functions.
\newblock In J.~J. Benedetto, editor, {\em Modern Sampling Theory}, chapter~7,
  pages 155--171. Birkh\"{o}user Boston, 2001.

\bibitem{GrochenigStrohmerMarvasti}
K.~Gr\"{o}chenig and T.~Strohmer.
\newblock Numerical and theoretical aspects of non-uniform sampling of
  band-limited images.
\newblock In F.~Marvasti, editor, {\em Nonuniform Sampling: Theory and
  Applications}, chapter~6, pages 283--324. Kluwer {A}cademic, {D}ordrecht,
  {T}he {N}etherlands, 2001.

\bibitem{PruessmannUnserMRIFast}
M.~Guerquin-Kern, M.~H{\"a}berlin, K.~P. Pruessmann, and M.~Unser.
\newblock {A fast wavelet-based reconstruction method for {M}agnetic
  {R}esonance {I}maging}.
\newblock {\em IEEE Trans. Med. Imaging}, 30(9):1649--1660, 2011.

\bibitem{HealyWeaverWaveletIEEE}
D.~M. Healy and J.~B. Weaver.
\newblock Two applications of wavelet transforms in {M}agnetic {R}esonance
  {I}maging.
\newblock {\em IEEE Trans. Inform. Theory}, 38(2):840--862, 1992.

\bibitem{hrycakIPRM}
T.~Hrycak and K.~Gr\"{o}chenig.
\newblock Pseudospectral {F}ourier reconstruction with the modified inverse
  polynomial reconstruction method.
\newblock {\em J. Comput. Phys.}, 229(3):933--946, 2010.

\bibitem{JacksonEtAlGridding}
J.~I. Jackson, C.~H. Meyer, D.~G. Nishimura, and A.~Macovski.
\newblock {Selection of a convolution function for {F}ourier inversion using
  gridding}.
\newblock {\em IEEE Trans. Med. Imaging}, 10:473--478, 1991.

\bibitem{Jaffard}
S.~Jaffard.
\newblock {A density criterion for frames of complex exponentials}.
\newblock {\em Mich. Math. J.}, 38(3):339--348, 1991.

\bibitem{shizgalGegen2}
J.-H. Jung and B.~D. Shizgal.
\newblock Generalization of the inverse polynomial reconstruction method in the
  resolution of the {G}ibbs phenomenon.
\newblock {\em J. Comput. Appl. Math.}, 172(1):131--151, 2004.

\bibitem{NMRbook}
J.~Keeler.
\newblock {\em Understanding NMR Spectroscopy}.
\newblock Wiley, 2010.

\bibitem{StanfordMRI}
A.~B. Kerr, J.~M. Pauly, B.~S. Hu, K.~C. Li, C.~J. Hardy, C.~H. Meyer,
  A.~Macovski, and D.~G. Nishimura.
\newblock {Real-time interactive {MRI} on a conventional scanner}.
\newblock {\em Magn. Reson. Med.}, 38(3):355--367, 1997.

\bibitem{KnoppKunisPotts}
T.~Knopp, S.~Kunis, and D.~Potts.
\newblock A note on the iterative {MRI} reconstruction from nonuniform k-space
  data.
\newblock {\em Int. J. Bio. Imag.}, page 24727, 2007.

\bibitem{KunisPottsScattered}
S.~Kunis and D.~Potts.
\newblock Stability results for scattered data interpolation by trigonometric
  polynomials.
\newblock {\em {SIAM} J. Sci. Comput.}, 29:1403--1419, 2007.

\bibitem{LaineWaveletsBiomed}
A.~F. Laine.
\newblock {Wavelets in temporal and spatial processing of biomedical images}.
\newblock {\em Annu. Rev. Biomed. Eng.}, 02:511--550, 2000.

\bibitem{ProlatesIII}
H.~J. Landau and H.~O. Pollak.
\newblock Prolate spheroidal wave functions, {F}ourier analysis and
  uncertainty---{III}: The dimension of the space of essentially time‐and
  band‐limited signals.
\newblock {\em Bell System Tech J.}, 41(4):1295--1336, 1962.

\bibitem{NMR2}
C.~Ludwig and U.~Gunther.
\newblock Metabolab - advanced nmr data processing and analysis for
  metabolomics.
\newblock {\em BMC Bioinformatics}, 12(1):366, 2011.

\bibitem{mallat09wavelet}
S.~G. Mallat.
\newblock {\em {A Wavelet Tour of Signal Processing: The Sparse Way}}.
\newblock Academic Press, 3 edition, 2009.

\bibitem{MarzoMZ}
J.~Marzo.
\newblock {Marcinkiewicz--Zygmund inequalities and interpolation by spherical
  harmonics}.
\newblock {\em J. Funct. Anal.}, 250(2):559--587, 2007.

\bibitem{MatejFesslerIterative}
S.~Matej, J.~A. Fessler, and I.~G. Kazantsev.
\newblock Iterative tomographic image reconstruction using {F}ourier-based
  forward and back-projectors.
\newblock {\em IEEE Trans. Med. Imaging}, 23:401--412, 2004.

\bibitem{NMR1}
P.~Mercier, M.~J. Lewis, D.~Chang, D.~Baker, and D.~S. Wishart.
\newblock Towards automatic metabolomic profiling of high-resolution
  one-dimensional proton nmr spectra.
\newblock {\em Journal of Biomolecular NMR}, 49(3-4):307--323, 2011.

\bibitem{MeyerEtAlSpiralCoronary}
C.~H. Meyer, B.~S. Hu, D.~G. Nishimura, and A.~Macovski.
\newblock {Fast spiral coronary artery imaging}.
\newblock {\em Magn. Reson. Med.}, 28:202--213, 1992.

\bibitem{NowakWaveletDenoise}
R.~Nowak.
\newblock {Wavelet-based {R}ician noise removal for {M}agnetic {R}esonance
  {I}maging}.
\newblock {\em IEEE Trans. Image Proc.}, 8:1408--19, 1998.

\bibitem{OrtegaCerdaMZ}
J.~Ortega-Cerda and J.~Saludes.
\newblock {Marcinkiewicz--Zygmund inequalities}.
\newblock {\em J. Approx. Theory}, 145:237--252, 2007.

\bibitem{PanychWaveletEncoding}
L.~P. Panych.
\newblock Theoretical comparison of {F}ourier and wavelet encoding in
  {M}agnetic {R}esonance {I}maging.
\newblock {\em IEEE Trans. Med. Imaging}, 15(2):141--153, 1996.

\bibitem{PanychEtAlWaveletEncoded}
L.~P. Panych, P.~D. Jakab, and F.~A. Jolesz.
\newblock Implementation of wavelet-encoded {MR} imaging.
\newblock {\em J. Magn. Reson. Imaging}, 3:649--55, 1993.

\bibitem{PottsEtAlNUFFTtutorial}
D.~Potts, G.~Steidl, and M.~Tasche.
\newblock Fast {F}ourier {T}ransforms for nonequispaced data: a tutorial.
\newblock In J.~J. Benedetto and P.~Ferreira, editors, {\em Modern Sampling
  Theory: Mathematics and Applications}, chapter~12, pages 249--274.
  Birkh{\"a}user, 2001.

\bibitem{Potts}
D.~Potts and M.~Tasche.
\newblock {Numerical stability of nonequispaced fast Fourier transforms}.
\newblock {\em J. Comput. Appl. Math.}, 222(2):655--674, 2008.

\bibitem{RosenfeldURS1}
D.~Rosenfeld.
\newblock {An optimal and efficient new gridding algorithm using singular value
  decomposition}.
\newblock {\em Magn. Reson. Med.}, 40(1):14--23, 1998.

\bibitem{SedaratDCFOptimal}
H.~Sedarat and D.~G. Nishimura.
\newblock {On the optimality of the gridding reconstruction algorithm}.
\newblock {\em IEEE Trans. Med. Imaging}, 19(4):306--317, 2000.

\bibitem{SeipJFA}
K.~Seip.
\newblock {On the connection between exponential bases and certain related
  sequences in {$L^2(-\pi,\pi)$}}.
\newblock {\em J. Funct. Anal.}, 130:131--160, 1995.

\bibitem{StrangNguyen}
G.~Strang and T.~Nguyen.
\newblock {\em Wavelets and Filter Banks}.
\newblock Wellesley-Cambridge Press, Wellesley, MA, 1996.

\bibitem{SunAvgFRI}
Q.~Sun.
\newblock Non-uniform average sampling and reconstruction of signals with
  finite rate of innovation.
\newblock {\em SIAM J. Math. Anal.}, 38:1389--1422, 2006.

\bibitem{SunXhouAverages}
W.~Sun and X.~Zhou.
\newblock Reconstruction of band-limited functions from local averages.
\newblock {\em Constr. Approx.}, 18:205--222, 2002.

\bibitem{FesslerFastIterativeMRI}
B.~P. Sutton, D.~C. Noll, and J.~A. Fessler.
\newblock {Fast, iterative image reconstruction for MRI in the presence of
  field inhomogeneities}.
\newblock {\em IEEE Trans. Med. Imaging}, 22(2):178--188, 2003.

\bibitem{unser1994general}
M.~Unser and A.~Aldroubi.
\newblock A general sampling theory for nonideal acquisition devices.
\newblock {\em IEEE Trans. Signal Process.}, 42(11):2915--2925, 1994.

\bibitem{GelbNonuniformFourier}
A.~Viswanathan, A.~Gelb, D.~Cochran, and R.~Renaut.
\newblock {On reconstructions from non-uniform spectral data}.
\newblock {\em J. Sci. Comput.}, 45(1--3):487--513, 2010.

\bibitem{WeaverEtAlWaveletFiltering}
J.~B. Weaver, Y.~Xu, D.~M. Healy, and J.~R. Driscoll.
\newblock {Filtering {MR} images in the wavelet transform domain}.
\newblock {\em Magn. Reson. Med.}, 21:288--295, 1991.

\bibitem{WeaverEtAlWaveletEncoding}
J.~B. Weaver, Y.~Xu, D.~M. Healy, and J.~R. Driscoll.
\newblock Wavelet-encoded {MR} imaging.
\newblock {\em Magn. Reson. Med.}, 24:275--287, 1992.

\bibitem{MRI1D3}
A.~B. Wolbarst, P.~Capasso, and A.~R. Wyant.
\newblock {\em MRI in One Dimension and with No Relaxation: A Gentle
  Introduction to a Challenging Subject}, pages 307--351.
\newblock John Wiley \& Sons, Inc., 2013.

\bibitem{young}
R.~M. Young.
\newblock {\em {An Introduction to Nonharmonic Fourier Series}}.
\newblock Academic Press Inc., San Diego, CA, first edition, 2001.

\end{thebibliography}

\end{document}